\documentclass[preprint]{elsarticle}

\usepackage{lineno,hyperref}
\modulolinenumbers[5]

\usepackage{amsmath}
\usepackage{amsfonts}
\usepackage{empheq}
\usepackage{stmaryrd}
\usepackage{physics}
\usepackage{mathtools}
\usepackage{subcaption}
\usepackage[export]{adjustbox}
\usepackage{yfonts}
\usepackage{mathrsfs}
\usepackage{url}

\journal{\url{www.arXiv.org}}

\begin{document}

\begin{frontmatter}

\title{Discontinuous Galerkin methods for a dispersive wave hydro-sediment-morphodynamic model}

\author[Oden]{Kazbek Kazhyken\corref{corrauth}}
\cortext[corrauth]{Corresponding author}
\ead{kazbek@oden.utexas.edu}

\author[CAMGSD]{Juha Videman}
\author[Oden]{Clint Dawson}

\address[Oden]{Oden Institute for Computational Engineering and Sciences, The University of Texas at Austin, Austin, TX 78712, USA}
\address[CAMGSD]{CAMGSD/Departamento de Matem\'atica, Universidade de Lisboa, Universidade de Lisboa, 1049-001 Lisbon, Portugal}

\begin{abstract}
A dispersive wave hydro-sediment-morphodynamic model developed by complementing the shallow water hydro-sediment-morphodynamic (SHSM) equations with the dispersive term from the Green-Naghdi equations is presented. A numerical solution algorithm for the model based on the second-order Strang operator splitting is presented. The model is partitioned into two parts, (1) the SHSM equations and (2) the dispersive correction part, which are discretized using discontinuous Galerkin finite element methods. This splitting technique provides a facility to  select dynamically regions of a problem domain where the dispersive term is not applied, e.g. wave breaking regions where the dispersive wave model is no longer valid. Algorithms that can handle wetting-drying and detect wave breaking are provided and a number of numerical examples are presented to validate the developed numerical solution algorithm. The results of the simulations indicate that the model is capable of predicting sediment transport and bed morphodynamic processes correctly provided that the empirical models for the suspended and bed load transport are properly calibrated. Moreover, the developed model is able to accurately capture hydrodynamics and wave dispersion effects up to swash zones, and its application is justified for simulations where dispersive wave effects are prevalent.
\end{abstract}


\end{frontmatter}


\section{Introduction}

A sediment transport process in coastal applications is a type of a two-phase fluid-solid flow with sea water as the fluid and pebbles and stones of varying sizes, and quartz sand as the solid. There are three modes of sediment transport: bed load, suspended load, and wash load transport. The bed load transport is characterized by motion of the sediment particles without detaching from the sediment bed for a significant amount of time, i.e. the sediment particles move by sliding, rolling, and saltating. There are a number of empirical models developed for the bed load transport, for example Meyer-Peter and Mueller \cite{meyer_and_muller_1948}, Fernandez Luque and Van Beek \cite{luque_beek_1976}, Nielsen \cite{nielsen_1992}, Ribberink \cite{ribberink_1998}. In the suspended load transport, the sediment particles suspended in water are advected with the water flow. These sediment particles, which are typically of a fine silt and clay size, remain suspended in water by turbulent flows and require a significant amount of time to settle on the sediment bed. Sediment particles in the wash load are transported without deposition while remaining close to the water surface in near-permanent suspension. Due to a limited effect of the wash load on the sediment bed morphology, effects of the wash load transport are not considered in the presented work.

Hydrodynamic, sediment transport, and bed morphodynamic processes are closely interrelated: hydrodynamic parameters of a water flow affect sediment transport rates, these rates influence the bed morphology that in its turn affects the water flow and sediment transport. These hydro-sediment-morphodynamic processes driven by astronomical tides, winds, and long-wave currents in coastal areas attract a high degree of interest since morphological changes of a coastal area can negatively affect its infrastructure and environment. Elements of coastal infrastructure, such as bridges, piers, and levees, can become structurally compromised as a result of excessive erosion of the sediment bed due to scouring. Environmental concerns include shoreline and beach erosion that may damage natural habitats of endangered protected species, and the effect of sediment transport on contaminants, i.e. sediment deposits may serve as dangerous contaminant sinks or sources. It is thus evident that mathematical modeling of hydro-sediment-morphodynamic processes in coastal areas has clear engineering relevance. Deriving such models poses, however, a number of challenges since they have to couple non-linear hydrodynamic, sediment transport, and bed morphodynamic equations along with modeling their two-way interactions. 

A number of hydro-sediment-morphodynamic models, ranging from one to three dimensional models, have been developed for coastal applications over the last four decades. These models are discussed in detail in \cite{amoudry_2008} and \cite{amoudry_souza_2011}. A three-dimensional model has the capacity for a more accurate and detailed resolution of the process \cite{wu_etal_2000, fang_and_wang_2000, marsooli_and_wu_2015}; however, the amount of computational resources required to run any sizable simulation with such a model is prohibitively large. Therefore, application of three-dimensional models is typically limited to short-time simulations over small-size domains. As an alternative, a depth averaged two- or, in some cases, one-dimensional model can be used to resolve hydro-sediment-morphodynamic processes in coastal areas. One such model is formed by the shallow water hydro-sediment-morphodynamic (SHSM) equations, which are derived by integrating and averaging the three-dimensional mass and momentum conservation equations of motion (e.g. see Wu \cite{wu_2007}). In the SHSM equations, the nonlinear shallow water equations, which resolve water-sediment mixture hydrodynamics, are fully coupled with sediment transport and bed morphodynamic models (see Cao \emph{et al.} \cite{cao_etal_2017} for variations of the SHSM equations). Within the last decade, the SHSM equations have been successfully applied in studies of coastal hydro-sediment-morphodynamic processes (e.g. Xiao \emph{et al.}, 2010 \cite{xiao_etal_2010}, Zhu and Dodd, 2015 \cite{zhu_and_dodd_2015}, Kim, 2015 \cite{kim_2015}, Incelli \emph{et al.}, 2016 \cite{incelli_etal_2016}, Briganti \emph{et al.}, 2016 \cite{briganti_etal_2016}). 

Numerical solution algorithms for the SHSM equations are typically developed with finite volume methods for applications with unstructured grids. Cao \emph{et al.}  \cite{cao_etal_2004} use the total-variation-diminishing (TVD) weighted average flux method (WAF) in conjunction with the Harten-Lax-van Leer-contact (HLLC) approximate Riemann solver to develop their numerical solution algorithm for the SHSM equations. Examples of works that employ HLLC as an approximate Riemann solver for numerical flux definitions include \cite{zhao_etal_2016}, \cite{zhao_etal_2019}, and \cite{hu_etal_2019}. Algorithms based on upwinding numerical fluxes and Roe-averaged states are developed in \cite{li_and_duffy_2011} and \cite{benkhaldoun_etal_2013}. Liu \emph{et al.} \cite{liu_etal_2015_0}, \cite{liu_etal_2015_1}, \cite{liu_and_beljadid_2017} develop numerical methods for the SHSM equations that employ a central-upwind scheme along with the Lagrange theorem to approximate the upper and lower bounds of the local wave speeds. Xia \emph{et al.}  \cite{xia_etal_2017} use the operator-splitting technique for the source term and the FORCE (first-order centered) approximate Riemann solver for a numerical treatment of the model. Discontinuous Galerkin discretizations of the SHSM equations are used less often, see, \emph{e.g.}, \cite{kesserwani_etal_2014} and \cite{clare_etal_2020}. 

The nonlinear shallow water equations, which form the hydrodynamic part of the SHSM equations, have a number of advantages: a capacity to approximate water motion with a sufficient accuracy in the shallow water flow regime, a plethora of developed numerical solution algorithms (e.g. Zhao \emph{et al.} \cite{zhao_etal_1994}, Anastasiou and Chan \cite{anastasiou_and_chan_1997}, Sleigh \emph{et al.}, \cite{sleigh_etal_1998}, Aizinger and Dawson \cite{aizinger_and_dawson_2002}, Yoon and Kang \cite{yoon_and_kang_2004},  Kubatko \emph{et al.} \cite{kubatko_etal_2006_nswe}, Samii \emph{et al.} \cite{samii_etal_2019}), efficient parallelization  strategies (e.g. hybrid MPI+OpenMP and HPX parallelization in Bremer \emph{et al.} \cite{bremer_etal_2019}), and its ability to approximate wave breaking effects in surf zones. However, this hydrodynamic model does not have a capacity to capture wave dispersion effects; and, therefore, an application of the SHSM equations is not feasible in areas where the dispersion effects are prevalent. An alternative depth-averaged hydrodynamic model that can reproduce dispersion effects is formed by the Green-Naghdi equations developed in \cite{green_naghdi_1976}. A number of numerical solution algorithms exist for the Green-Naghdi equations that use various discretization techniques, from finite difference to finite element methods, and a Strang operator splitting technique (e.g. see \cite{chazel_etal_2011, bonneton_etal_2011, panda_etal_2014, lannes_and_marche_2015, duran_marche_2015, duran_and_marche_2017, samii_and_dawson_2018, marche_2020}). The use of a Strang operator splitting in these algorithms provides the capacity to switch between the nonlinear shallow water equations and the Green-Naghdi equations whenever one of the hydrodynamic models is more accurate than the other \cite{duran_marche_2015}.

The purpose of the presented work is to introduce dispersive wave effects into the SHSM equations. This is achieved through considering the Green-Naghdi equations  which results in a dispersive wave hydro-sediment-morphodynamic model. Since the difference between the nonlinear shallow water equations and the Green-Naghdi equations is constituted by the dispersive term defined through a differential operator that forms an elliptic system \cite{bonneton_etal_2011}, this new model is formed by incorporating the dispersive term into the SHSM equations. The resulting model has the potential to be used in the simulation of morphodynamic processes in areas where dispersive wave effects are prevalent. Numerical solution algorithms for this model are developed employing a Strang operator splitting technique and discontinuous Galerkin finite element methods. A significant portion of this work comprises the development of a massively parallel solver that uses the developed numerical solution algorithms.  The solver extends a C++ software package developed by Bremer and Kazhyken\footnote{The software is under development on the date of the publication, and can be accessed at \url{www.github.com/UT-CHG/dgswemv2}. Should there be any questions, comments, or suggestions, please contact the developers through the repository issues page.}. 

The rest of the paper is organized as follows. Section 2 presents the governing equations for the dispersive wave hydro-sediment-morphodynamic model. The developed numerical solution algorithms are introduced in Section 3. Section 4 presents a number of numerical tests, including one-dimensional and two-dimensional dam break simulations and solitary wave runs over an erodible sloping beach, that are used to perform verification and validation of the developed algorithms. Final conclusions are presented in Section 5. 

\section{Governing equations}
\begin{figure}
\center
\includegraphics[width=3in]{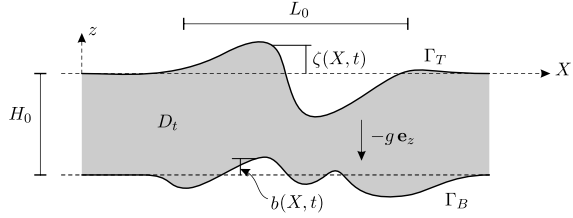}
\caption{A model representation of a body of water as a domain $D_t \subset \mathbb R^{d+1}$.}
\label{Fig:Domain}
\end{figure}

A body of water can be represented by a domain $D_t \subset \mathbb R^{d+1}$, where $d$ is the horizontal spatial dimension that can take values 1 or 2, and $t$ represents the time variable. The domain $D_t$ is filled with a water-sediment mixture, modeled as an incompressible inviscid fluid, and bounded vertically by the bottom and top boundaries, $\Gamma_B$ and $\Gamma_T$, which the fluid particles cannot cross (cf. Fig.\ref{Fig:Domain}). It is assumed that $\Gamma_B$ and $\Gamma_T$ can be represented as graphs that vary in time: $\Gamma_B$  due to sediment transport and bed morphodynamic processes, $\Gamma_T$  as the evolving free surface of the body of water. The bathymetry, $b(X,t)$, and the free surface elevation, $\zeta(X,t)$, of the body of water are used in the parameterization of $\Gamma_B$ and $\Gamma_T$: 
\begin{linenomath}
\begin{subequations}
\begin{align}
\Gamma_B &= \{(X,-H_0+b(X,t)):X\in \mathbb R^d\}, \\
\Gamma_T &= \{(X,\zeta(X,t)):X\in \mathbb R^d\},
\end{align}
\end{subequations}
\end{linenomath}
and the domain $D_t$ is defined as a set of points $(X,z) \in \mathbb R^d \times \mathbb R$ where $-H_0+b(X,t) < z < \zeta(X,t)$.

A depth-averaged model that can resolve water wave dynamics, and subsequent sediment transport and bed evolution in the domain $D_t$ is the shallow water hydro-sediment-morphodynamic (SHSM) equations (e.g. see Cao \emph{et al.} \cite{cao_etal_2004}). The hydrodynamic part of the equations is represented by the nonlinear shallow water equations, which provide a sufficiently accurate approximation to the water wave dynamics whenever the shallowness parameter $\mu=H_0^2/L_0^2$, where $L_0$ is the characteristic length, and $H_0$ is the reference depth, is less than unity. The present work aims to develop a hydro-sediment-morphodynamic model that has the capacity to capture wave dispersion effects, which the nonlinear shallow water equations are unable to resolve. Therefore, the nonlinear shallow water equations in the SHSM model are replaced with a single parameter variation of the Green-Naghdi equations, a depth-averaged hydrodynamic model which has the capacity to capture wave dispersion effects, introduced by Bonneton \emph{et al.} in \cite{bonneton_etal_2011}. This forms a set of equations defined over a horizontal domain $\Omega \subset \mathbb R^d$:
\begin{linenomath}
\begin{equation}\label{Eq:GNHSM}
\partial_t \boldsymbol q + \nabla \cdot \boldsymbol F(\boldsymbol q) + \boldsymbol{D}(\boldsymbol q)  = \boldsymbol S(\boldsymbol q),
\end{equation}
\end{linenomath}
where the vector of unknowns $\boldsymbol q$ and the flux matrix $\boldsymbol F(\boldsymbol q)$ are
\begin{linenomath}
\begin{equation}
\boldsymbol q = \begin{Bmatrix} h \\ h \mathbf u \\ hc \\ b \end{Bmatrix}, \quad
\boldsymbol F(\boldsymbol q) = \begin{Bmatrix} h \mathbf u \\ h\mathbf u \otimes \mathbf u + \frac 1 2 g h^2 \mathbf I \\ h c \mathbf u \\ \mathbf q_b \end{Bmatrix},
\end{equation}
\end{linenomath}
the source term $\boldsymbol S(\boldsymbol q)$ is defined as
\begin{linenomath}
\begin{equation}
\boldsymbol S(\boldsymbol q) = \begin{Bmatrix}  \frac{E-D}{1-p}  \\ -gh \nabla b - \frac{\rho_s-\rho_w}{2\rho}gh^2 \nabla c-\frac{(\rho_0-\rho)(E-D)}{\rho(1-p)} \mathbf u + \mathbf f \\ E-D \\  - \frac{E-D}{1-p} \end{Bmatrix},
\end{equation}
\end{linenomath}
$\mathbf u$ is the water velocity represented by a $d$ dimensional vector and $h$ is the water depth represented by the mapping $h(X,t) = \zeta(X,t) + H_0 - b(X,t)$ and assumed to be bounded from below by a positive value. Moreover,  $c$ is the volume concentration of sediment in water-sediment mixture, $E$ and $D$ are the sediment entrainment and deposition rates, respectively, $p$ is the bed porosity, $\rho_w$ and $\rho_s$ are the water and the sediment densities, $\rho$ and $\rho_0$ are the water-sediment mixture and saturated bed densities defined as $\rho=(1-c)\rho_w+c\rho_s$ and $\rho_0=(1-p)\rho_s+p\rho_w$, $\mathbf q_b$ is the bed load sediment flux, $\mathbf f$ comprises additional source terms for the  momentum continuity equation (e.g. the Coriolis, bottom friction, and surface wind stress forces), $g$ is the acceleration due to gravity, and $\mathbf I \in \mathbb R^{d\times d}$ is the identity matrix. Finally, the wave dispersion effects are introduced into the model through the dispersive term 
\begin{linenomath}
\begin{equation}
\boldsymbol D(\boldsymbol q) = \begin{Bmatrix} 0 \\ \mathbf w_1 - \alpha^{-1} g h \nabla \zeta \\ 0 \\ 0 \end{Bmatrix},
\end{equation}
\end{linenomath}
where $\mathbf w_1$ is defined through an elliptic system
\begin{linenomath}
\begin{equation}\label{Eq:w1}
(\mathbf I + \alpha h \mathcal T h^{-1}) \mathbf w_1 = \alpha^{-1} g h \nabla \zeta + h \mathcal Q_1(\mathbf u),
\end{equation}
\end{linenomath}
with operators $\mathcal T$ and $\mathcal Q_1$ defined as
\begin{linenomath}
\begin{subequations}
\begin{align}
\mathcal T(\mathbf w) =&\,\mathcal R_1(\nabla\cdot\mathbf w) + \mathcal R_2(\nabla b \cdot \mathbf w), \\
\mathcal Q_1(\mathbf w) =&-2\mathcal R_1\left(\partial_{x} \mathbf w \cdot \partial_{y} \mathbf w^\perp+(\nabla \cdot \mathbf w)^2\right)+\mathcal R_2\left(\mathbf w\cdot (\mathbf w \cdot \nabla)\nabla b\right),
\end{align}
\end{subequations}
\end{linenomath}
where operators $\mathcal R_1$ and $\mathcal R_2$ are
\begin{linenomath}
\begin{subequations}
\begin{align}
\mathcal R_1(w) &= -\frac 1 {3h} \nabla(h^3 w) -  \frac {h}{2} w \nabla b,\\
\mathcal R_2(w) &= \frac 1 {2h} \nabla (h^2 w) +  w \nabla b,
\end{align}
\end{subequations}
\end{linenomath}
and $\mathbf w^\perp = (-w_2, w_1)^\mathbf{T}$. Parameter $\alpha\in\mathbb R$ in the dispersive term is used to optimize dispersive properties of the presented hydro-sediment-morphodynamic model. By adjusting $\alpha$, the difference between the phase and group velocities coming from the Stokes linear theory and the Green-Naghdi equations can be minimized. A common strategy aims at minimizing the averaged variation over some range of wave number values \cite{bonneton_etal_2011}.

In the presented model $E$, $D$ and $\mathbf q_b$ are defined through empirical equations. The sediment entrainment rate $E$ may be defined as in \cite{li_duffy_2011}:
\begin{linenomath}
\begin{equation}
E=
\begin{cases}
\phi (\theta-\theta_c)\lvert\mathbf{u}\rvert h&\text{if}\,\,\,\, \theta>\theta_c \\
0&\text{if}\,\,\,\,\theta\leq\theta_c 
\end{cases},
\end{equation}
\end{linenomath}
where $\phi$ is a calibration parameter, $\theta_c$ is the critical Shields parameter and $\theta$ is the Shields parameter given by $\theta=\lvert\boldsymbol{\tau}_b\rvert/\sqrt{sgd_{50}}$, where $\boldsymbol{\tau}_b$ is the bottom friction, $s=\rho_s/\rho_w-1$ is the submerged specific gravity, and $d_{50}$ is the mean sediment particle size. The sediment deposition rate $D$ can be estimated by an empirical model from \cite{cao_etal_2004}:
\begin{linenomath}
\begin{equation}
D = \omega_o C_a(1-C_a)^2,
\end{equation}
\end{linenomath}
where $\omega_o$ is the setting velocity of a sediment particle in still water, and $C_a = c\alpha_c$ is the near-bed sediment volume concentration with the coefficient $\alpha_c = \min(2, (1-p)/c)$. A number of empirical models for $\mathbf q_b$  is proposed as (see \cite{diaz_etal_2008, cordier_etal_2011} and all the references therein)
\begin{linenomath}
\begin{equation}\label{Eq:Qb}
\mathbf q_b = A(h, \mathbf u)\mathbf u \lvert \mathbf u \rvert^{m-1},
\end{equation}
\end{linenomath}
where $1 \leq m \leq 3$ and $A(h, \mathbf u)$ is an empirical equation, e.g. the Grass model takes $A$ as a constant calibrated for the application under investigation and sets $m=3$, cf. \cite{grass_1981}.

\section{Numerical methods}

Discontinuous Galerkin finite element methods are used to discretize the governing equations. This choice facilitates the use of unstructured meshes that are well suited for irregular geometries of coastal areas. Thus, the problem domain $\Omega$ is partitioned into a finite element mesh $\mathcal T_h = \{K\}$ that provides an approximation to the domain: 
\begin{linenomath}
\begin{equation}
\Omega\approx\Omega_h=\sum_{K\in \mathcal{T}_h}K,
\end{equation}
\end{linenomath}
where the subscript $h$ stands for the mesh parameter represented by the diameter of the smallest element in the mesh. The set of all mesh element faces, $\partial\mathcal T_h$, and the set of all edges of the mesh skeleton, $\mathcal{E}_h$, are defined as
\begin{linenomath}
\begin{subequations}
\begin{align}
\partial\mathcal T_h &= \lbrace\partial K : K\in \mathcal{T}_h\rbrace,\\
\mathcal{E}_h &= \lbrace e\in\bigcup_{K\in\mathcal{T}_h}\partial K\rbrace.
\end{align}
\end{subequations}
\end{linenomath}
Note that in $\mathcal E_h$ the common element faces appear only once but in $\partial \mathcal T_h$ they are counted twice. 

To develop variational formulations of the governing equations, inner products are defined for finite dimensional vectors $\boldsymbol u$ and $\boldsymbol v$ through:
\begin{linenomath}
\begin{subequations}
\begin{align}
(\boldsymbol u,\boldsymbol v)_\Omega &= \int_\Omega \boldsymbol u \cdot \boldsymbol v \, \dd X, \\
\langle \boldsymbol u, \boldsymbol v \rangle_{\partial \Omega} &= \int_{\partial \Omega}\boldsymbol u \cdot \boldsymbol v\, \dd X,
\end{align}
\end{subequations}
\end{linenomath}
for $\Omega \subset \mathbb R^d$ and $\partial \Omega \subset \mathbb R^{d-1}$.

An approximating space of trial and test functions is chosen as the set of square integrable functions over $\Omega_h$ such that their restriction to an element $K$ belongs to $\mathcal Q^p(K)$, a space of polynomials of degree at most $p \ge 0$ with support in $K$:
\begin{linenomath}
\begin{equation}
\mathbf V_h^{p,m} \coloneqq \{\boldsymbol v \in (L^2(\Omega_h))^{m}: \boldsymbol v|_K \in \mathcal (\mathcal Q^p(K))^{m} \quad  \forall K \in \mathcal T_h \},
\end{equation}
\end{linenomath}
and, similarly, an approximation space over the mesh skeleton is chosen as
\begin{linenomath}
\begin{equation}
\mathbf M_h^{p,m} \coloneqq \{\boldsymbol \mu \in (L^2(\mathcal E_h))^{m}: \boldsymbol \mu|_e \in \mathcal (\mathcal Q^p(e))^{m} \quad \forall e \in \mathcal E_h\}.
\end{equation}
\end{linenomath}

A Strang operator splitting technique is used in the numerical solution of the hydro-sediment-morphodynamic model presented in Eq.(\ref{Eq:GNHSM}). To this end, the model is split into two separate parts: (1) the SHSM equations obtained by dropping the dispersive term of the equations, and (2) the dispersive correction part where the wave dispersion effects on flow velocities are introduced into the model through the dispersive term. If $\mathcal S_1$ is a numerical solution operator for the SHSM equations, i.e. $\mathcal S_1(\Delta t)$ propagates numerical solution by a time step $\Delta t$, and, similarly, $\mathcal S_2$ is a numerical solution operator for the dispersive correction part, then the numerical solution operator for the full hydro-sediment-morphodynamic model in Eq.(\ref{Eq:GNHSM}) can be approximated with the Strang operator splitting technique \cite{strang_1968}:
\begin{linenomath}
\begin{equation}
\mathcal S(\Delta t) = \mathcal S_1(\Delta t/2) \mathcal S_2(\Delta t) \mathcal S_1(\Delta t/2),
\end{equation}
\end{linenomath}
where $\mathcal S$ is a second-order temporal discretization if both $\mathcal S_1$ and $\mathcal S_2$ use a second-order time discretization method.

A numerical solution operator $\mathcal S_1$ for the SHSM equations is developed using a discontinuous Galerkin finite element formulation where an approximate solution $\boldsymbol q_h \in \mathbf V_h^{p,d+3}$ is sought such that it satisfies the variational formulation
\begin{linenomath}
\begin{equation}\label{Eq:NSWEVarLoc}
(\partial_t \boldsymbol{q}_h,\boldsymbol{v})_{\mathcal{T}_h}-(\boldsymbol{F}_h,\nabla \boldsymbol{v})_{\mathcal{T}_h}+\langle\boldsymbol{F}_h^*,\boldsymbol{v}\rangle_{\partial {\mathcal{T}_h}}-(\boldsymbol{S}_h,\boldsymbol{v})_{\mathcal{T}_h}=0\quad \forall \boldsymbol{v}\in\textbf{V}_h^{p,d+3},
\end{equation}
\end{linenomath}
where $\boldsymbol{F}_h=\boldsymbol F(\boldsymbol{q}_h)$ and $\boldsymbol{S}_h=\boldsymbol S(\boldsymbol{q}_h)$, ${\boldsymbol F_h^*}$ is a single valued approximation to $\boldsymbol F_h \mathbf n$ over element faces, called the numerical flux, and $\mathbf n$ is the unit outward normal vector to element face. To define the numerical flux, the bed update part of the SHSM equations is singled out for a separate treatment. The numerical flux for this formulation is then defined as 
\begin{equation}
\boldsymbol{F}_h^*=\begin{Bmatrix}\boldsymbol{G}_h^* \\ \mathbf{q}_b^*\end{Bmatrix},
\end{equation}
where $\mathbf{q}_b^*$ is the numerical bed load flux, and $\boldsymbol{G}_h^*$ is the numerical flux for the remaining part of the system where the vector of unknowns $\boldsymbol r$ and the flux matrix $\boldsymbol G(\boldsymbol r)$ are
\begin{linenomath}
\begin{equation}
\boldsymbol r = \begin{Bmatrix} h \\ h \mathbf u \\ hc \end{Bmatrix}, \quad
\boldsymbol G(\boldsymbol r) = \begin{Bmatrix} h \mathbf u \\ h\mathbf u \otimes \mathbf u + \frac 1 2 g h^2 \mathbf I \\ h c \mathbf u  \end{Bmatrix}.
\end{equation}
\end{linenomath}

Assuming that the sediment transport is always in the flow direction, the numerical flux $\mathbf q_b^*$ is defined as in \cite{mirabito_etal_2011}:
\begin{linenomath}
\begin{equation}
\mathbf q_b^*=
\begin{cases}
\mathbf q_b^+&\text{if}\,\,\,\,\mathbf {\hat u} \cdot \mathbf n \geq 0 \\
\mathbf q_b^-&\text{if}\,\,\,\,\mathbf {\hat u} \cdot \mathbf n < 0 
\end{cases},
\end{equation}
\end{linenomath}
where $\mathbf {\hat u}$ is the Roe-averaged velocity defined as
\begin{linenomath}
\begin{equation}
\mathbf {\hat u} = \frac{\mathbf{u}^+\sqrt{h^+} + \mathbf{u}^-\sqrt{h^-}}{\sqrt{h^+}+\sqrt{h^-}}.
\end{equation}
\end{linenomath}
Here and for the rest of this article, superscript $+$ denotes a variable value at $\partial K$ when approaching from the interior of an element $K$, and $-$ when approaching from the exterior. An upwinding scheme is employed for the numerical bed load flux $\mathbf q_b^*$ since computing the eigenvalues of the normal Jacobian matrix for the flux matrix $\boldsymbol F(\boldsymbol{q})$ requires computationally intensive numerical approximation techniques and does not guarantee real values except in the case where the Grass model is used for $\mathbf q_b$ \cite{diaz_etal_2008, cordier_etal_2011}. Therefore, using numerical flux definitions that involve the eigenvalues of the normal Jacobian matrix for the full system may prove to be unfeasible.

The normal Jacobian matrix $\boldsymbol A = \partial_{\boldsymbol r} (\boldsymbol G \mathbf n)$ of the remaining part of the system has four real eigenvalues: $\lambda_{1,2} = \mathbf{u}\cdot\mathbf{n}\pm\sqrt{gh}$, $\lambda_{3,4} = \mathbf{u}\cdot\mathbf{n}$. A Godunov-type Harten–Lax–van Leer scheme is used to define the numerical flux for the remaining system \cite{harten_etal_1983}:
\begin{linenomath}
\begin{equation}
\boldsymbol{G}_h^* = \begin{cases}
\boldsymbol{G}_h^+\mathbf{n}&\text{if}\,\,\,\,S^+>0\\
\boldsymbol{G}_h^{\text{HLL}}&\text{if}\,\,\,\,S^+\leq0\leq S^-\\
\boldsymbol{G}_h^-\mathbf{n}&\text{if}\,\,\,\,S^-<0
\end{cases},
\end{equation}
\end{linenomath}
where $\boldsymbol{G}_h=\boldsymbol G(\boldsymbol{r}_h)$, the truncated characteristic speeds $S^+$ and $S^-$ are 
\begin{linenomath}
\begin{subequations}
\begin{align}
S^+&=\min(\mathbf{u}^+\cdot\mathbf{n}-\sqrt{gh^+}, \mathbf{u}^-\cdot\mathbf{n}-\sqrt{gh^-}),\\
S^-&=\max(\mathbf{u}^+\cdot\mathbf{n}+\sqrt{gh^+}, \mathbf{u}^-\cdot\mathbf{n}+\sqrt{gh^-}),
\end{align}
\end{subequations}
\end{linenomath}
and the Harten–Lax–van Leer flux $\boldsymbol{G}_h^{\text{HLL}}$ is
\begin{linenomath}
\begin{equation}
\boldsymbol{G}_h^{\text{HLL}}=\frac{1}{S^--S^+}((S^-\boldsymbol{G}_h^+-S^+\boldsymbol{G}_h^-)\mathbf{n}-S^+S^-(\boldsymbol{r}_h^+-\boldsymbol{r}_h^-)).
\end{equation}
\end{linenomath}

A hybridized discontinuous Galerkin scheme may be used to define the numerical flux through $\widehat{\boldsymbol r}_h\in \mathbf M_h^{p,d+2}$, an approximation to ${\boldsymbol r}$ over the mesh skeleton called the numerical trace \cite{nguyen_peraire_2012}:
\begin{linenomath}
\begin{equation}
\boldsymbol G_h^* = \widehat{{\boldsymbol G}}_h\mathbf n + \boldsymbol \tau ({\boldsymbol r}_h - \widehat{\boldsymbol r}_h),
\end{equation}
\end{linenomath}
where $\widehat{\boldsymbol G}_h=\boldsymbol G(\widehat{\boldsymbol r}_h)$, and $\boldsymbol \tau =\lambda_{\max}(\widehat{\boldsymbol r}_h)$ is the stabilization parameter defined as the maximum eigenvalue of the normal Jacobian matrix $\boldsymbol A$:
\begin{linenomath}
\begin{equation}
\lambda_{\max}(\boldsymbol r) = \lvert\mathbf{u}\cdot\mathbf{n}\rvert+\sqrt{gh}.
\end{equation}
\end{linenomath} 
The numerical trace $\widehat {\boldsymbol r}_h\in \mathbf M_h^{p,d+2}$ must be such that the numerical flux is conserved across all internal edges in the mesh skeleton, and boundary conditions are satisfied at all boundary edges through the boundary operator $\boldsymbol B_h$ defined according to an imposed boundary condition \cite{nguyen_peraire_2012}: 
\begin{linenomath}
\begin{equation}\label{Eq:NSWEVarGlob}
\langle \boldsymbol G^*_h , \boldsymbol \mu \rangle_{\partial \mathcal T_h \backslash \partial \Omega_h} + \langle \boldsymbol B_h , \boldsymbol \mu \rangle_{\partial \mathcal T_h \cap \partial \Omega_h}=0\,\,\,\forall \boldsymbol\mu\in \mathbf M_h^{p,d+2}.
\end{equation}
\end{linenomath}
Eq.(\ref{Eq:NSWEVarLoc}) and Eq.(\ref{Eq:NSWEVarGlob}) along with the definition of $\mathbf q_b^*$ form a system of equations that is used to solve for an approximate solution $\boldsymbol q_h \in \mathbf V_h^{p,d+3}$. The boundary condition operator $\boldsymbol B_h$ is defined as \begin{linenomath}
\begin{equation}
\boldsymbol B_h = \boldsymbol A^+ \boldsymbol r_h - |\boldsymbol A| \widehat{\boldsymbol r}_h - \boldsymbol A^- {\boldsymbol r}_\infty,
\end{equation}
\end{linenomath}
where $\boldsymbol A^{\pm} = \frac{1}{2}(\boldsymbol A \pm |\boldsymbol A|)$, and ${\boldsymbol r}_\infty$ is the weakly imposed boundary state \cite{nguyen_peraire_2012}
. For a slip wall boundary condition, $\boldsymbol B_h$ is defined as
\begin{linenomath}
\begin{equation}
\boldsymbol B_h =  \widehat{\boldsymbol r}_h - {\boldsymbol r}_{\text{slip}},
\end{equation}
\end{linenomath}
where ${\boldsymbol r}_{\text{slip}} = \{(h)_h \quad (h\mathbf u)_h - ((h\mathbf u)_h \cdot \textbf{n}) \textbf{n}  \quad (hc)_h\}^\mathbf{T}$ is a state with its normal velocity component truncated \cite{nguyen_peraire_2012}.

In order to generate $\mathcal S_2$, a numerical solution operator for the dispersive correction part of the presented hydro-sediment-morphodynamic model, Eq.(\ref{Eq:w1}) is written as a system of first order equations using the definition for operator $\mathcal T$ \cite{samii_and_dawson_2018}:
\begin{linenomath}
\begin{empheq}[left=\empheqlbrace]{equation}\label{Eq:w1w2}
\begin{split}
&\nabla \cdot (h^{-1} \mathbf w_1) - h^{-3} w_2 = 0\\
&\mathbf w_1- \tfrac 1 3 \nabla w_2 - \tfrac {1}{2} h^{-1} w_2 \nabla b 
+ \tfrac 1 2 \nabla (h \nabla b \cdot \mathbf w_1) +
 \mathbf w_1 \nabla b \otimes \nabla b = \mathbf{s}(\boldsymbol{q})
\end{split},
\end{empheq}
\end{linenomath}
where $\textbf{s}(\boldsymbol{q}) = \alpha^{-1}gh \nabla \zeta + h \mathcal Q_1(\mathbf u)$. A discontinuous Galerkin finite element discretization for Eq.(\ref{Eq:w1w2}) forms a global system of equations. A hybridized discontinuous Galerkin formulation can be used to reduce the dimension of the global system of equations. Therefore, the hybridized discontinuous Galerkin method developed by Samii and Dawson in  \cite{samii_and_dawson_2018} is employed to treat numerically Eq.(\ref{Eq:w1w2}) to obtain an approximate solution $\mathbf{w}_{1h}\in\textbf{V}_h^{p,d}$. The result is then used in the dispersive correction to seek an approximate solution $\boldsymbol{q}_h\in\textbf{V}_h^{p,d+3}$ that satisfies the variational formulation
\begin{linenomath}
\begin{equation}
(\partial_t \boldsymbol{q}_h,\boldsymbol{v})_{\mathcal{T}_h}+\left(\boldsymbol D_h,\boldsymbol{v}\right)_{\mathcal{T}_h}=0\quad \forall \boldsymbol{v}\in\textbf{V}_h^{p,d+3},
\end{equation}
\end{linenomath}
where $\boldsymbol D_h = \boldsymbol D(\boldsymbol{q}_h)$. High order derivatives of $\mathbf{u}_h$, present in $\mathcal Q_1(\mathbf{u}_h)$, are computed weakly using a discontinuous Galerkin method with centered numerical fluxes. 

In the developed depth-averaged hydro-sediment-morphodynamic model, it is assumed that the water depth $h$ is bounded from below by a positive value. This assumption implemented by a wetting-drying algorithm which ensures that the water depth remains positive. The numerical solution operator $\mathcal S_2$ does not affect the water depth; therefore, the wetting-drying algorithm should work in conjunction with the numerical solution operator for the SHSM equations $\mathcal S_1$. In the presented work, the wetting-drying algorithm developed for the nonlinear shallow water equations by Bunya \emph{et al.} in \cite{bunya_etal_2009} is adapted to the SHSM equations. In the adapted version of the Bunya \emph{et al.} wetting-drying algorithm, the sediment term $hc$ in the SHSM equations is treated the same way as the momentum term $h\mathbf{u}$  and the rest of the algorithm remains the same. The bed update part of the equations does not affect the water depth and, therefore, it does not require the wetting-drying algorithm. Finally, in the dispersive correction part of the equations the wet-dry front is modeled as a slip wall boundary.

Using the Green-Naghdi equations as the hydrodynamic part of the presented model allows capturing wave dispersion effects; however, the Green-Naghdi equations are limited to parts of the problem domain that are free from discontinuities in numerical solutions \cite{duran_and_marche_2017}. This poses certain limitations on the application of the Green-Naghdi equations, e.g. wave breaking phenomena in surf zones present themselves as a water depth discontinuity in numerical solutions. While the Green-Naghdi equations cannot accurately resolve wave breaking, the nonlinear shallow water equations are more suitable for such areas \cite{duran_and_marche_2017}. Using the Strang operator splitting allows switching to the nonlinear shallow water equations from the Green-Naghdi equations by setting $\mathcal{S}_2=1$ in regions with discontinuities in numerical solutions. Thus, a discontinuity detection criterion is required to dynamically switch to $\mathcal{S}_2=1$. In the presented work, the numerical solution algorithm is augmented with the water depth discontinuity detection criterion adopted by Duran and Marche in \cite{duran_and_marche_2017} from Krivodonova \emph{et al.} \cite{krivodonova_etal_2004}. A water depth discontinuity is identified over an element $K$ if the parameter \cite{krivodonova_etal_2004,duran_and_marche_2017}
\begin{linenomath}
\begin{equation}
\mathbb I_K = \frac{\sum_{F\in\partial K_{\text{in}}}\vert\int_{F}(h^+-h^-)\dd X\vert}{\mathfrak{h}_K^{\frac{p+1}{2}}\,\vert \partial K_{\text{in}}\vert\,\Vert h\Vert_{L^{\infty}(K) }}
\end{equation}
\end{linenomath}
is greater than a specified threshold that is typically $O(1)$. In this description of the parameter $\mathbb I_K$, $\mathfrak{h}_K$ is the element diameter,  $\partial K_{\text{in}}$ are the inflow faces of the element where $\mathbf{u}\cdot\mathbf{n}<0$, and $\vert \partial K_{\text{in}}\vert$ is the total length of the inflow faces. 

Since $\mathcal{S}_2$ is not applied in regions with discontinuities in  the numerical solutions, a slope limiting is not needed for the dispersive correction part of the presented model. However, whenever discontinuities occur in the numerical solutions to the SHSM equations a slope limiting algorithm is required in order to remove the oscillations at sharp discontinuities and to preserve numerical stability. Thus, the Cockburn-Shu limiter \cite{cockburn_shu_2001} is incorporated into the numerical solution algorithm and applied in conjunction with the operator $\mathcal S_1$. The details of the limiter are not presented here, but readers are encouraged to consult the original source.    

\section{Numerical experiments and discussion}

The developed numerical model has been implemented in a software framework written in C++ programming language with the use of open source scientific computing libraries, such as Eigen \cite{eigen}, Blaze \cite{blaze}, and PETSc \cite{petsc}. The software has been parallelized for shared and distributed memory systems with the use of a hybrid OpenMP+MPI programming, and HPX \cite{hpx}. Performance comparison between the hybrid programming and HPX has been performed by Bremer \emph{et al.} in \cite{bremer_etal_2019}. 

The presented numerical model is validated in five numerical examples. In the first four set-up only the numerical solution operator for the SHSM equations is validated against four dam break experiments. In these experiments the dispersive wave effects are negligible; therefore, $\mathcal S_2=1$ in the simulations. The last example uses the full dispersive wave hydro-sediment-morphodynamic model to simulate water waves, sediment transport, and bed morphodynamics caused by solitary wave runs over a sloping beach.

The first-order Dubiner polynomials from \cite{dubiner_1991} are used for the approximating space $\mathbf V_h$, and the first-order Legendre polynomials are used for the approximating space $\mathbf M_h$. In all presented  examples, numerical solutions are computed using two different definitions of the numerical flux $\boldsymbol{G}^*_h$: (1) the Harten–Lax–van Leer discontinuous Galerkin scheme (HLL DG), (2) the Nguyen-Peraire hybridized discontinuous Galerkin scheme (NP HDG). Consequently, the numerical results obtained using these two definitions for the numerical flux are compared against each other.

\subsection{1D dam break}
\begin{figure}[!t]
\center
\includegraphics[trim=0.75in 0.25in 0.75in 0.25in,clip,width=4.75in]{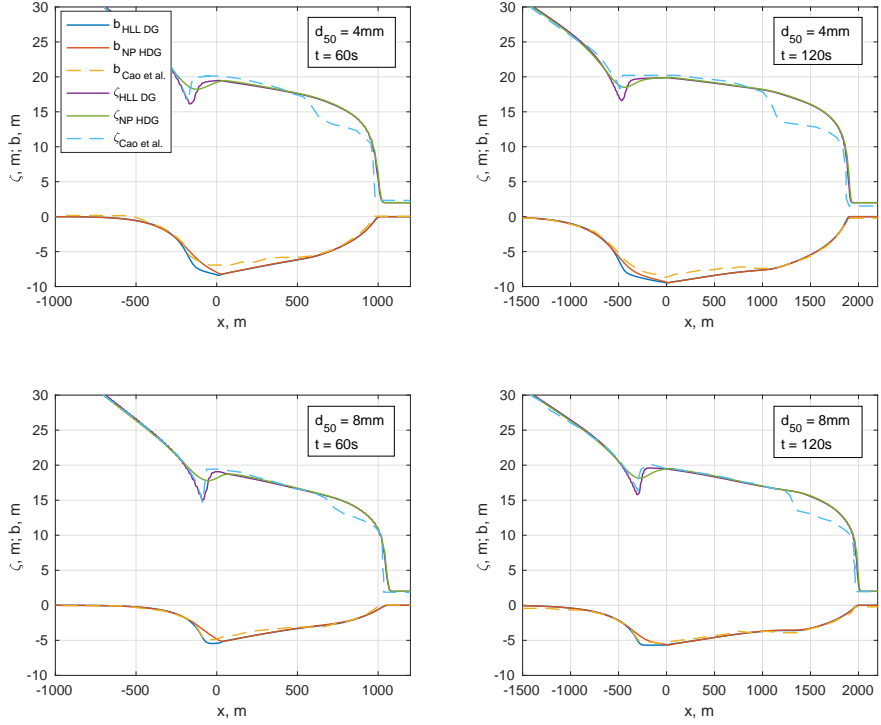}
\caption{Free surface elevation and bathymetry from the 1D dam break simulation compared with Cao \emph{et al.} experiment \cite{cao_etal_2004}.}
\label{Fig:EX1}
\end{figure}

In this numerical experiment the SHSM equations are used to simulate a 1D dam break over a mobile bed. Initial conditions for this experiment are set as a clear ($c_0(x)=0$) still water ($\mathbf{u}_0(x)=0$) with its depth distributed as
\begin{linenomath}
\begin{equation}
h_0(x) = \begin{cases}
40&\text{if}\,\,\,\,x\leq0\\
2&\text{if}\,\,\,\,x>0
\end{cases},
\end{equation}
\end{linenomath}
and the bathymetry set to $b_0(x)=0$. The mobile bed in this experiment has the sediment density $\rho_s=2650\,\text{kg}/\text{m}^3$, the bed porosity $p=0.4$, the critical Shields parameter $\theta_c=0.045$, and the mean sediment particle size $d_{50}$ set as 4mm and 8mm for two separate simulation runs. For the sediment entrainment rate model, the calibration parameter is set as $\phi=0.015$. The bed load transport is not considered in this numerical experiment by setting $\mathbf{q}_b=0$. The bottom friction force is introduced into the model through the source term $\boldsymbol S(\boldsymbol q)$ by setting
\begin{linenomath}
\begin{equation}
\mathbf f = \frac{gn^2}{h^{1/3}}{\vert\mathbf{u}\vert\mathbf{u}},
\label{Eq:f}
\end{equation}
\end{linenomath}
with the Manning’s roughness coefficient $n=0.03$.

The problem domain $\Omega = (-5000,5000)\times(-10,10)\,\text{m}^2$ is partitioned into a finite element mesh with $500\times1$ square cells each containing 2 triangular elements. The explicit Euler time stepping scheme is employed with the time step $\Delta t = 0.1$s. Two simulations with varying mean sediment particle sizes are run for 2 minutes, and their results are compared to the numerical experiments carried out for the same 1D dam break problem by Cao \emph{et al.} in \cite{cao_etal_2004}. The results of the numerical simulations at $t=\{60, 120\}$s for $d_{50}=\{4, 8\}$mm are presented in Fig.\ref{Fig:EX1}. Smaller sediment particle sizes imply larger magnitude for sediment entrainment rate $E$, which presents itself as a larger bed erosion for $d_{50}=4\text{mm}$. The numerical results for both the free surface elevation, $\zeta$, and the bathymetry, $b$, are in good agreement with the results obtained by Cao \emph{et al.} The numerical results obtained with HLL DG and NP HDG schemes closely match each other except in the area of the hydraulic jump where NP HDG scheme provides a smoother solution for the free surface elevation. 

\subsection{1D dam break with wetting-drying}
\begin{figure}[!t]
\center
\includegraphics[trim=0.75in 0.5in 0.75in 0.5in,clip,width=4.75in]{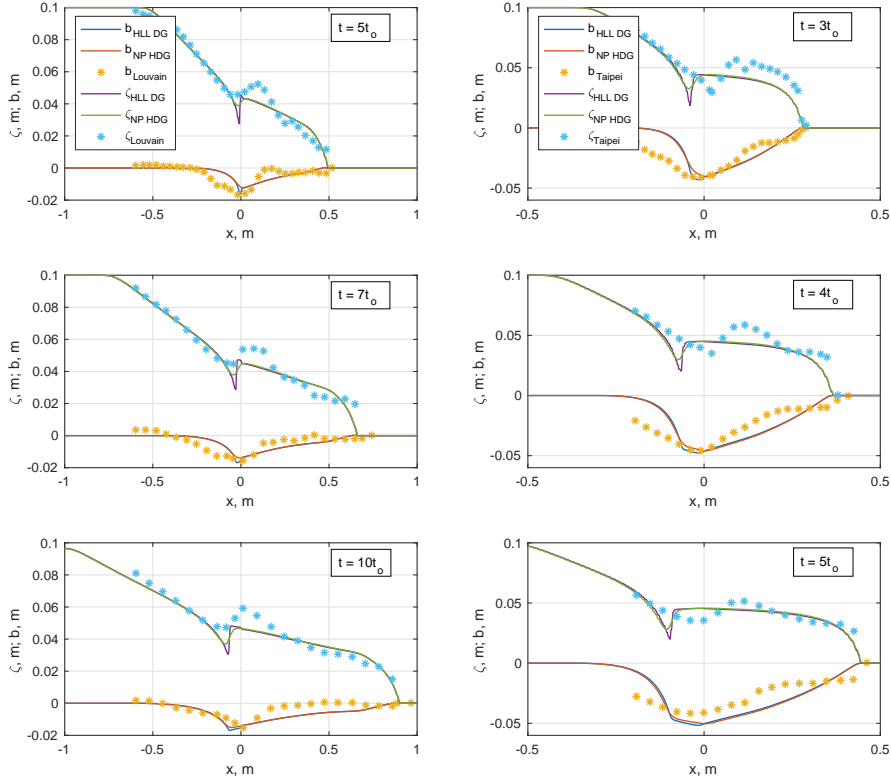}
\caption{Free surface elevation and bathymetry from the 1D dam break with wetting-drying simulations compared with the Louvain \cite{fraccarollo_capart_2002} and Taipei \cite{capart_young_1998} experiments.}
\label{Fig:EX2}
\end{figure}

This example simulates a 1D dam break over a mobile dry bed and is used to validate the wetting-drying algorithm employed in the presented numerical model. Numerical simulations for this experiment are performed with the SHSM equations where water is initially in clear still state, the water depth is set to
\begin{linenomath}
\begin{equation}
h_0(x) = \begin{cases}
0.1&\text{if}\,\,\,\,x\leq0\\
0&\text{if}\,\,\,\,x>0
\end{cases},
\end{equation}
\end{linenomath}
and the initial bathymetry is $b_0(x)=0$. Two physical experiments have been performed for this setup: (1) the Louvain experiment by Fraccarollo and Capart \cite{fraccarollo_capart_2002}, (2) the Taipei experiment by Capart and Young \cite{capart_young_1998}. These experiments are  set up similarly except for the sediment properties. In the Louvain experiment the sediment density $\rho_s=1540\,\text{kg}/\text{m}^3$, the bed porosity $p=0.3$, the critical Shields parameter $\theta_c=0.05$, and the mean sediment particle size $d_{50}=3.5$mm. On the other hand, in the Taipei experiment the sediment density $\rho_s=1048\,\text{kg}/\text{m}^3$, the bed porosity $p=0.28$, the critical Shields parameter $\theta_c=0.05$, and the mean sediment particle size $d_{50}=6.1$mm. The calibration parameter for the sediment entrainment rate model, $\phi$, is set as 4.0 for the Louvain experiment, and 2.5 for the Taipei experiment. In both experiments, the bed load transport is disregarded by setting $\mathbf{q}_b=0$, and the Manning's friction model from Eq.(\ref{Eq:f}) is used for the bottom friction force with $n=0.025$.

The problem domain $\Omega = (-1,1)\times(-2\cdot10^{-3},2\cdot10^{-3})\,\text{m}^2$ is partitioned into a finite element mesh with $500\times1$ square cells each containing two triangular elements. The explicit Euler time stepping scheme with the time step $\Delta t = 5\cdot10^{-4}$s is used to propagate simulations in time for 1s. The simulations of the 1D dam break over mobile dry bed are carried out with the parameters from the Louvain and Taipei experiments. The results are compared with the Louvain experiment at $t=\{5t_0,7t_0,10t_0\}$ and with the Taipei experiment at $t=\{3t_0,4t_0,5t_0\}$, where $t_0=\sqrt{g/h_0}\approx0.101$s ($h_0=0.1$m), in Fig.\ref{Fig:EX2}. The numerical solution algorithm successfully models the wetting-drying process while providing sufficiently accurate numerical results for the free surface elevation, $\zeta$, and the bathymetry, $b$. Similar to the previous example, HLL DG and NP HDG results closely match each other everywhere other than the hydraulic jump area.

\subsection{2D flume with abrupt widening}
\begin{figure}[t!]
\center
\includegraphics[trim=0.25in 0.25in 0.25in 0.25in,clip,width=4.75in]{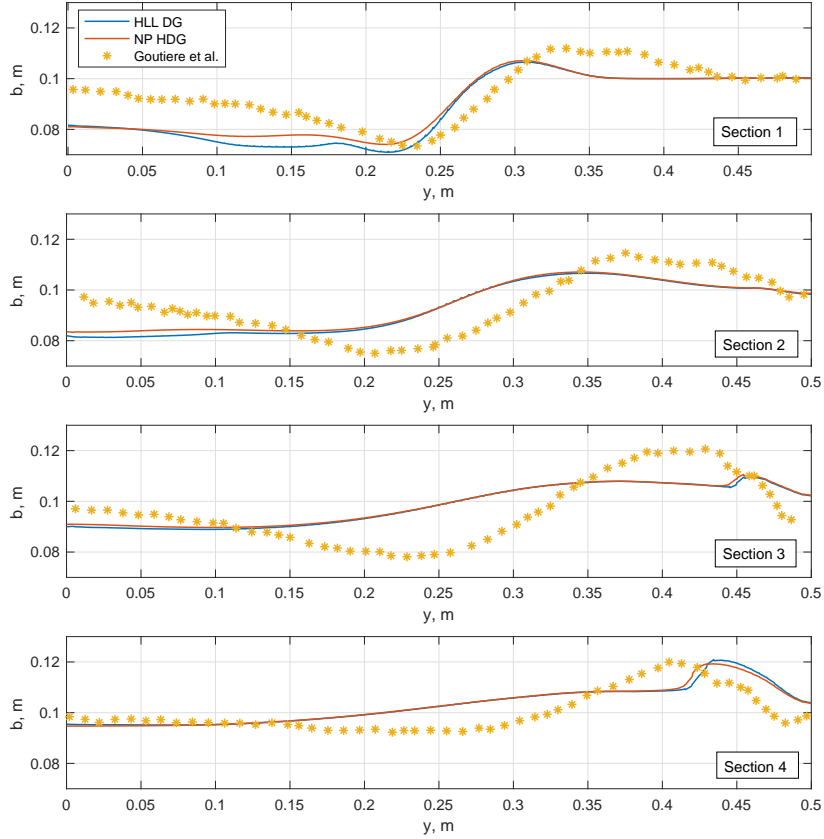}
\caption{Sediment erosion/deposition measurements from the 2D flume with abrupt widening experiment compared with Goutiere \emph{et al.} results \cite{goutiere_etal_2011}.}
\label{Fig:EX3}
\end{figure}

A 2D dam break is simulated in an "L-shaped" flume which is 0.25m wide in its initial 4m and has an abrupt widening on one side to 0.5m for the remaining 2m. The flume bed is covered with 0.1m of sediment ($b_0(x)=0.1$) with the following properties: the sediment density $\rho_s=2630\,\text{kg}/\text{m}^3$, the bed porosity $p=0.39$, the critical Shields parameter $\theta_c=0.047$, the mean sediment particle size $d_{50}=1.72$mm. In this experiment, only the suspended load in taken into account while setting the calibration parameter for the sediment entrainment rate model, $\phi$, to 0.35. Initial conditions for the SHSM equations simulations are clear still water with its initial depth 
\begin{linenomath}
\begin{equation}
h_0(x) = \begin{cases}
0.25&\text{if}\,\,\,\,x\leq3\\
0&\text{if}\,\,\,\,x>3
\end{cases},
\end{equation}
\end{linenomath}
which implies that the abrupt expansion of the flume is located 1m downstream from the dam break location. The Manning's friction model from Eq.(\ref{Eq:f}) is used for the bottom friction force with $n=0.0165$.

The "L-shaped" problem domain $\Omega$ for this simulation is partitioned into nearly $4\cdot10^4$ triangular elements. The explicit Euler time integration scheme is used for this numerical simulation with the time step $\Delta t = 2\cdot10^{-4}$s. The simulation is run for 20s after which the sediment erosion/deposition measurements are taken at 4 lateral sections located at $x=\{4.1(\text{S1}), 4.2(\text{S2}), 4.3(\text{S3}), 4.4(\text{S4})\}$m. These measurements are compared with the results of the physical experiment performed by Goutiere \emph{et al.} in \cite{goutiere_etal_2011} in Fig.\ref{Fig:EX3}. The results of the numerical simulation generally agree with the results of the physical experiment. A general tendency for sediment erosion on the left side and sediment deposition on the right side of the flume is captured in the numerical simulation. The model is also able to capture large sediment deposition on the right side at Sections 3 and 4 where the water flow experiences sudden deceleration due to an impact with the side wall \cite{goutiere_etal_2011}. No significant differences can be observed between HLL DG and NP HDG schemes in this example.     

\subsection{2D partial dam break}
\begin{figure}[t!]
\center
\includegraphics[trim=0.25in 0.25in 0.25in 0.25in,clip,width=4.75in]{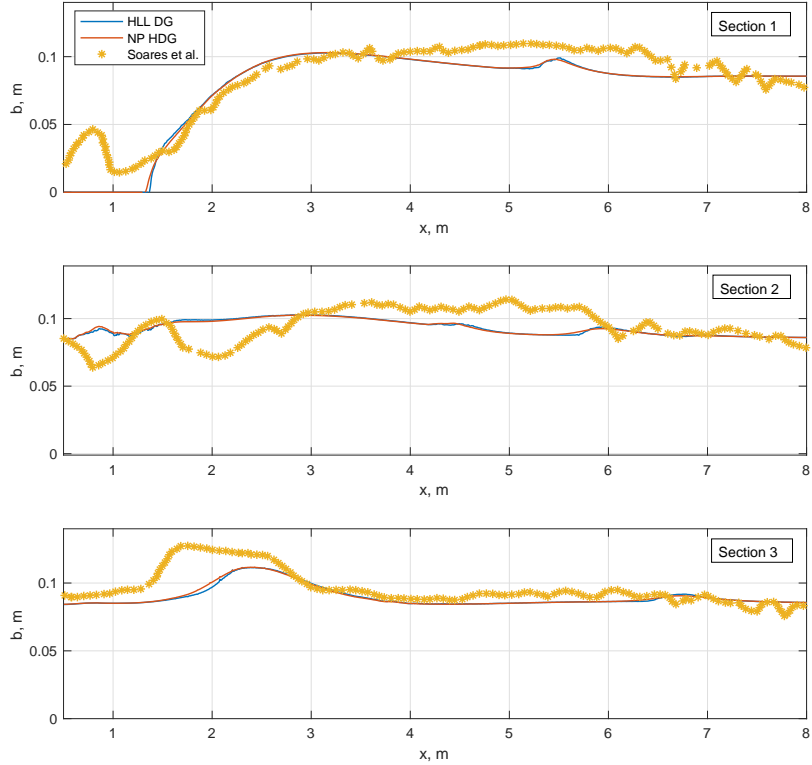}
\caption{Sediment erosion/deposition measurements from the 2D partial dam break experiment compared with Soares-Fraz\~ao \emph{et al.} results \cite{soares_etal_2012}.}
\label{Fig:EX4}
\end{figure}

A partial 2D dam break is simulated in a flume that consist of two 3.6m wide reservoirs that are connected by a 1m long and 1m wide channel with a gate in the middle, which is removed at the beginning of the experiment to simulate a partial dam break. The channel connects the reservoirs along their longitudinal axes. The wet reservoir that holds water is 10m long, and the dry reservoir is 15m long. The bed of the dry reservoir is covered by 0.085m of sediment with the sediment density $\rho_s=2630\,\text{kg}/\text{m}^3$, the bed porosity $p=0.42$, the critical Shields parameter $\theta_c=0.047$, and the mean sediment particle size $d_{50}=1.61$mm. The bed load transport is not taken into account in this experiment, and the calibration parameter for the sediment entrainment rate model $\phi=0.05$. Initially, the wet reservoir water is in clear still state and is 0.47m deep. The bottom friction force is modeled with the Manning's friction model from Eq.(\ref{Eq:f}) with $n=0.0165$. 

The problem domain $\Omega$ for this numerical experiment is partitioned into over $10^5$ triangular elements. The numerical simulation is propagated in time with the explicit Euler time stepping scheme with the time step $\Delta t = 5\cdot10^{-4}$s. After 20s of the numerical simulation, the sediment erosion/deposition measurements are taken at 3 longitudinal sections of the dry reservoir located at $y=\{0.2(\text{S1}), 0.7(\text{S2}), 1.45(\text{S3})\}$m away from the longitudinal axis of the reservoir. Fig.\ref{Fig:EX4} presents the measurements and compares them with the results of the physical experiment performed by Soares-Fraz\~ao \emph{et al.} in \cite{soares_etal_2012}. The results of the numerical simulation are in good agreement with the results of the physical experiment. The sediment is mostly eroded near the channel, where the bed is nearly completely scoured away and deposited downstream by the water flow from the dam break, as is evident from the measurements at Section 1. In this example, HLL DG and NP HDG schemes did not lead to significantly different numerical solutions.   

\subsection{Solitary wave over a sloping beach}
\begin{figure}[t!]
\center
\includegraphics[trim=0.5in 0.5in 0.5in 0.5in,clip,width=4.75in]{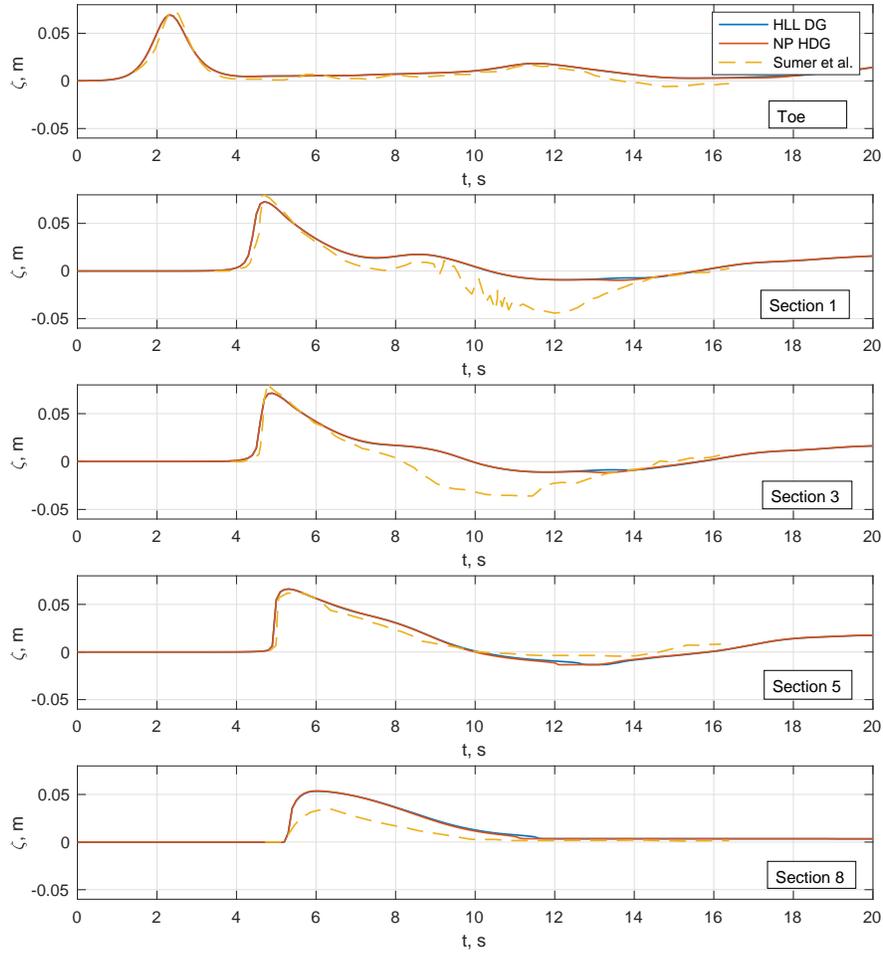}
\caption{Free surface elevation measurements at 5 measuring stations compared to the experimental results by Sumer \emph{et al.} \cite{sumer_etal_2011}.}
\label{Fig:Wave}
\end{figure}

In this experiment, the full dispersive wave hydro-sediment-morphodynamic model is used to simulate water waves, and subsequent sediment transport and bed evolution during run up and run down of a solitary wave over a linearly sloping beach. This experiment showcases a number of features of the presented model: (1) the use of the Green-Naghdi equations as a hydrodynamic component of the model since wave dispersion effects play a significant role during run up of a solitary wave over a sloping beach, (2) switching to the nonlinear shallow water equations as a hydrodynamic model in swash zones since solitary waves in this experiment have a sufficiently high amplitude to experience wave breaking, (3) solitary waves that run over a sloping beach in this experiment cause significant erosion/deposition of the beach bed; thus, the ability of the model to estimate sediment transport and bed morphology can be evaluated. Initial conditions for solitary waves in this experiment are characterized by equations 
\begin{linenomath}
\begin{equation}
h_0(x) = H_0 + a_0\sech^2\left(\kappa(x-x_0)\right), 
\quad
(h\mathbf u)_0(x) = c_0 h_0(x) - c_0 H_0,
\end{equation}
\end{linenomath}
where $a_0$ is the solitary wave height, $x_0$ the initial wave position, and
\begin{linenomath}
\begin{equation}
\kappa = \frac {\sqrt{3 a_0}} {2H_0 \sqrt{H_0+a_0}},
\quad
c_0 = \sqrt{g (H_0+a_0)}.
\end{equation}
\end{linenomath}

Initially, a simulation has been performed over a rigid bed to validate the dispersive wave hydrodynamic model. To carry out this numerical simulation, the problem domain $\Omega = (-10,10)\times(-2.5\cdot10^{-2},2.5\cdot10^{-2})\,\text{m}^2$ is partitioned into a finite element mesh comprised of $400\times1$ square cells containing two triangular elements. A two-stage second-order Runge-Kutta method is used to perform time integration with the time step $\Delta t=5\cdot10^{-3}$s. The Manning’s roughness coefficient $n=0.03$ is used for the bottom friction force. The toe of the sloping beach for this simulation is located at $x=0$ where an initially flat bed starts climbing linearly up at a 1:14 rate. The parameters for the solitary wave in this simulation are: $H_0=0.4$m, $a_0=0.071$m, and $x_0=-5$m. This simulation setup corresponds to the solitary wave run over a sloping beach experiment performed by Sumer \emph{et al.} \cite{sumer_etal_2011}. Fig.\ref{Fig:Wave} presents numerical solutions for the free surface elevations recorded at 5 measuring stations located at $x=\{0.0(\text{Toe}), 4.63(\text{S1}), 4.87(\text{S3}), 5.35(\text{S5}), 5.85(\text{S8})\}$m during 20s of the simulation and compares them to the experimental results provided by Sumer \emph{et al.} The  experimental results suggest that wave breaking occurs somewhere between Sections 3 and 5. This is accurately captured with the dispersive wave hydrodynamic model. However, the free surface elevation measurements at the onshore Section 8 show that the hydrodynamic model is less precise in resolving water waves in the swash zone. Subsequently, the hydrodynamic model is unable to simulate accurately the water motion during the run down stage. Nevertheless, considering complexities associated with modeling water motion induced by solitary waves over a sloping beach, the results of the simulation can be regarded as satisfactory. 

\begin{figure}[t!]
\center
\includegraphics[trim=0.5in 0.0in 0.5in 0.0in,clip,width=4.75in]{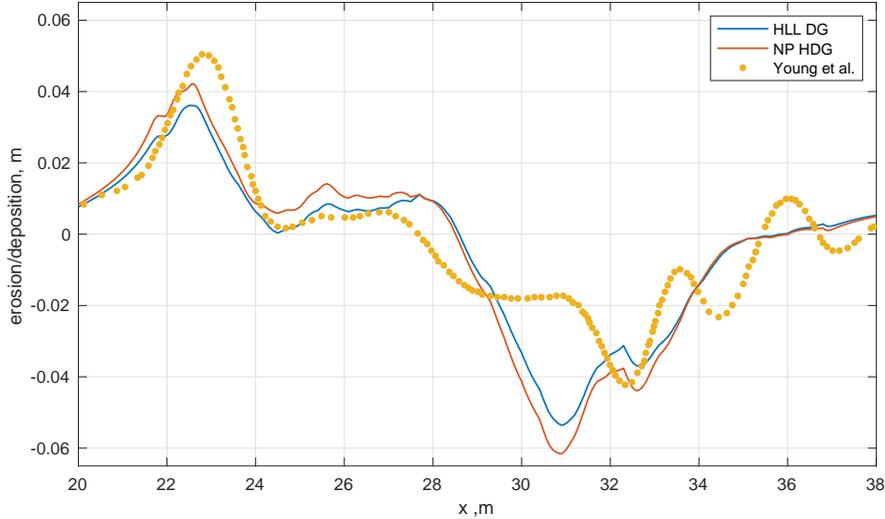}
\caption{Sediment erosion/deposition measurements for a simulation with the suspended load transport compared with the results by Young \emph{et al.} \cite{young_etal_2010}.}
\label{Fig:EX5_suspended}
\end{figure}

\begin{figure}[t!]
\center
\includegraphics[trim=0.5in 0.0in 0.5in 0.0in,clip,width=4.75in]{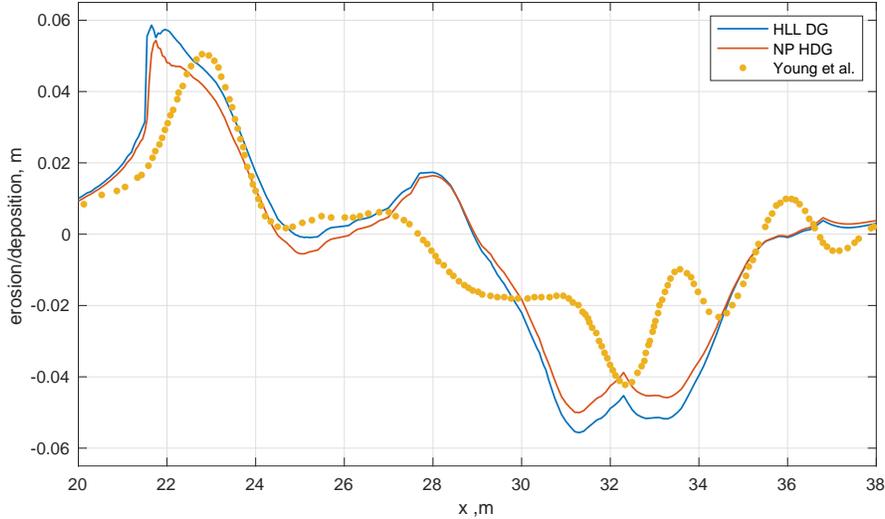}
\caption{Sediment erosion/deposition measurements for a simulation with the suspended and bed load transport compared with the results by Young \emph{et al.} \cite{young_etal_2010}.}
\label{Fig:EX5_full}
\end{figure}

To validate the sediment transport and bed morphodynamic part of the model, solitary wave run simulations have been performed over the problem domain $\Omega = (-8,42)\times(-5\cdot10^{-2},5\cdot10^{-2})\,\text{m}^2$. The problem domain is partitioned into $500\times1$ square cells each containing two triangular elements, and a two stage second-order Runge-Kutta method with the time step $\Delta t=2.5\cdot10^{-3}$s is used for temporal discretization. The toe of the sloping beach in the simulation is located at $x=12$m where the flat rigid bed starts climbing at 1:15 rate. The sloping part of the beach is covered with mobile sediment with the sediment density $\rho_s=2650\,\text{kg}/\text{m}^3$, the bed porosity $p=0.4$, the critical Shields parameter $\theta_c=0.045$, the mean sediment particle size $d_{50}=0.2$mm. The Manning’s roughness coefficient $n=0.008$ is used for the bottom friction force. The solitary wave in this simulation is parametrized with $H_0=1$m, $a_0=0.6$m, and $x_0=2$m. A physical experiment with the same setup has been performed by Young \emph{et al.} in \cite{young_etal_2010} where a number of solitary waves have been run over a sloping beach and subsequent sediment erosion/deposition has been recorded. Two simulations are performed: (1) a simulation where only the suspended load transport is taken into account with its results presented in Fig.\ref{Fig:EX5_suspended}, and (2) a simulation where both the suspended and bed load transport are considered with its results presented in Fig.\ref{Fig:EX5_full}. For the suspended load, the calibration parameter for the sediment entrainment rate model, $\phi$, is set to 0.35; and the Grass model with $A=2\cdot10^{-4}$ is used as a model for the bed load flux $\mathbf{q}_b$. In both of these simulations sediment erosion/deposition measurements are taken after 3 solitary waves have been run over the sloping beach for 2m each, which is a sufficient time for water to substantially settle. The results of these measurements are compared with the experimental results by Young \emph{et al.} and they are in good agreement. The experimental results indicate that \cite{young_etal_2010}: (1) during the initial run up sediment is entrained in water and deposited onshore at the maximum excursion point where the water flow stalls, (2) during the run down process a shallow high velocity flow causes net sediment erosion in the region between $x=24$m and $x=35$m, (3) this entrained sediment is then deposited offshore in the vicinity of the hydraulic jump, which is formed by the retreating water, due to sudden deceleration of the sediment-rich flow. The numerical model is able to capture the sediment transport and bed morphodynamics features observed in the experiment accurately.

\section{Conclusions}

A dispersive wave hydro-sediment-morphodynamic model has been developed by introducing the dispersive term of a single parameter variation of the Green-Naghdi equations into the SHSM equations. The model can be used to simulate water waves, and the resulting sediment transport and bed morphodynamic processes in areas where wave dispersion effects are prevalent. A numerical solution operator has been developed for the model which employs the second-order Strang operator splitting technique. In order to employ this technique, the dispersive term has been singled out for a separate numerical treatment with a hybridized discontinuous Galerkin method developed by Samii and Dawson in \cite{samii_and_dawson_2018}, and Harten–Lax–van Leer discontinuous Galerkin, and Nguyen-Peraire hybridized discontinuous Galerkin schemes have been developed for the remaining SHSM equations. The splitting technique makes it possible to select regions where the dispersive term is not applied, e.g. in wave breaking regions where the dispersive wave model is no longer valid. The numerical model is augmented with a wave breaking detection mechanism that can dynamically determine regions where the dispersive term is not applied. To facilitate the use of the developed model in problems where water may completely recede from parts of the problem domain, the wetting-drying algorithm by Bunya \emph{et al.} \cite{bunya_etal_2009} has been incorporated into the numerical model.

The numerical model has been validated against a number of numerical examples. Dam break simulations have been performed to validate the numerical solution schemes developed for the SHSM equations. The results of the simulations indicate that the developed schemes are able to capture hydro-sediment-morphodynamic processes with a sufficient accuracy. Since empirical models are used for the suspended and bed load transport, a close calibration for the empirical models' parameters may be required to improve the accuracy of the presented model. Simulations of a solitary wave run-up over a sloping beach have been performed to validate the full dispersive wave hydro-sediment-morphodynamic model. The results of the simulations indicate that the use of the presented model is justified for flows where the wave dispersion effects are prevalent. Subsequently, the use of the presented model for such flows accurately captures sediment transport and bed morphodynamic processes driven by these flows.     

\section{Acknowledgments}
This work has been supported by funding from the National Science Foundation Grant 1854986, and the Portuguese government through Funda\c{c}\~ao para a Ci\^encia e a Tecnologia (FCT), I.P., under the project DGCOAST (UTAP-EXPL/MAT/0017/2017). Authors would like to acknowledge the support of the Texas Advanced Computing Center through the allocation TG-DMS080016N used in the parallel computations of this work. 

\bibliography{mybibfile}

\begin{thebibliography}{10}
\expandafter\ifx\csname url\endcsname\relax
  \def\url#1{\texttt{#1}}\fi
\expandafter\ifx\csname urlprefix\endcsname\relax\def\urlprefix{URL }\fi
\expandafter\ifx\csname href\endcsname\relax
  \def\href#1#2{#2} \def\path#1{#1}\fi

\bibitem{meyer_and_muller_1948}
E.~Meyer-Peter, R.~M\"uller, Formulas for bed-load transport, Proceedings of
  2nd meeting of the International Association for Hydraulic Structures
  Research (1948) 39–64.

\bibitem{luque_beek_1976}
{R. Fernandez Luque}, {R. van Beek}, Erosion {A}nd {T}ransport {O}f
  {B}ed-{L}oad {S}ediment, Journal of Hydraulic Research 14~(2) (1976)
  127--144.
\newblock \href {http://dx.doi.org/10.1080/00221687609499677}
  {\path{doi:10.1080/00221687609499677}}.

\bibitem{nielsen_1992}
P.~Nielsen, Coastal Bottom Boundary Layers and Sediment Transport, Advanced
  series on ocean engineering, World Scientific, 1992.

\bibitem{ribberink_1998}
J.~S. Ribberink, Bed-load transport for steady flows and unsteady oscillatory
  flows, Coastal Engineering 34~(1) (1998) 59 -- 82.
\newblock \href {http://dx.doi.org/10.1016/S0378-3839(98)00013-1}
  {\path{doi:10.1016/S0378-3839(98)00013-1}}.

\bibitem{amoudry_2008}
L.~O. Amoudry, A {R}eview on {C}oastal {S}ediment {T}ransport {M}odelling,
  \url{http://nora.nerc.ac.uk/id/eprint/8360} (2008).

\bibitem{amoudry_souza_2011}
L.~O. Amoudry, A.~J. Souza, Deterministic {C}oastal {M}orphological and
  {S}ediment {T}ransport {M}odeling: a {R}eview and {D}iscussion, Reviews of
  Geophysics 49~(2).
\newblock \href {http://dx.doi.org/10.1029/2010RG000341}
  {\path{doi:10.1029/2010RG000341}}.

\bibitem{wu_etal_2000}
W.~Wu, W.~Rodi, T.~Wenka, {3D Numerical Modeling of Flow and Sediment Transport
  in Open Channels}, Journal of Hydraulic Engineering 126~(1) (2000) 4--15.
\newblock \href {http://dx.doi.org/10.1061/(ASCE)0733-9429(2000)126:1(4)}
  {\path{doi:10.1061/(ASCE)0733-9429(2000)126:1(4)}}.

\bibitem{fang_and_wang_2000}
H.-W. Fang, G.-Q. Wang, {Three-Dimensional Mathematical Model of
  Suspended-Sediment Transport}, Journal of Hydraulic Engineering 126~(8)
  (2000) 578--592.
\newblock \href {http://dx.doi.org/10.1061/(ASCE)0733-9429(2000)126:8(578)}
  {\path{doi:10.1061/(ASCE)0733-9429(2000)126:8(578)}}.

\bibitem{marsooli_and_wu_2015}
R.~Marsooli, W.~Wu, {Three-Dimensional Numerical Modeling of Dam-Break Flows
  with Sediment Transport over Movable Beds}, Journal of Hydraulic Engineering
  141~(1) (2015) 04014066.
\newblock \href {http://dx.doi.org/10.1061/(ASCE)HY.1943-7900.0000947}
  {\path{doi:10.1061/(ASCE)HY.1943-7900.0000947}}.

\bibitem{wu_2007}
W.~Wu, Computational River Dynamics, CRC Press, London, 2007.
\newblock \href {http://dx.doi.org/10.4324/9780203938485}
  {\path{doi:10.4324/9780203938485}}.

\bibitem{cao_etal_2017}
Z.~Cao, C.~Xia, G.~Pender, Q.~Liu, {Shallow Water Hydro-Sediment-Morphodynamic
  Equations for Fluvial Processes}, Journal of Hydraulic Engineering 143~(5)
  (2017) 02517001.
\newblock \href {http://dx.doi.org/10.1061/(ASCE)HY.1943-7900.0001281}
  {\path{doi:10.1061/(ASCE)HY.1943-7900.0001281}}.

\bibitem{xiao_etal_2010}
H.~Xiao, Y.~L. Young, J.~H. Prévost, {Hydro- and morpho-dynamic modeling of
  breaking solitary waves over a fine sand beach. Part II: Numerical
  simulation}, Marine Geology 269~(3) (2010) 119 -- 131.
\newblock \href {http://dx.doi.org/10.1016/j.margeo.2009.12.008}
  {\path{doi:10.1016/j.margeo.2009.12.008}}.

\bibitem{zhu_and_dodd_2015}
F.~Zhu, N.~Dodd, The morphodynamics of a swash event on an erodible beach,
  Journal of Fluid Mechanics 762 (2015) 110–140.
\newblock \href {http://dx.doi.org/10.1017/jfm.2014.610}
  {\path{doi:10.1017/jfm.2014.610}}.

\bibitem{kim_2015}
D.-H. Kim, {H2D} morphodynamic model considering wave, current and sediment
  interaction, Coastal Engineering 95 (2015) 20 -- 34.
\newblock \href {http://dx.doi.org/10.1016/j.coastaleng.2014.09.006}
  {\path{doi:10.1016/j.coastaleng.2014.09.006}}.

\bibitem{incelli_etal_2016}
G.~Incelli, N.~Dodd, C.~E. Blenkinsopp, F.~Zhu, R.~Briganti, Morphodynamical
  modelling of field-scale swash events, Coastal Engineering 115 (2016) 42 --
  57.
\newblock \href {http://dx.doi.org/10.1016/j.coastaleng.2015.09.006}
  {\path{doi:10.1016/j.coastaleng.2015.09.006}}.

\bibitem{briganti_etal_2016}
R.~Briganti, A.~Torres-Freyermuth, T.~E. Baldock, M.~Brocchini, N.~Dodd, T.-J.
  Hsu, Z.~Jiang, Y.~Kim, J.~C. Pintado-Patiño, M.~Postacchini, Advances in
  numerical modelling of swash zone dynamics, Coastal Engineering 115 (2016) 26
  -- 41.
\newblock \href {http://dx.doi.org/10.1016/j.coastaleng.2016.05.001}
  {\path{doi:10.1016/j.coastaleng.2016.05.001}}.

\bibitem{cao_etal_2004}
Z.~Cao, G.~Pender, S.~Wallis, P.~Prof, Computational dam-break hydraulics over
  erodible sediment bed, Journal of Hydraulic Engineering 130~(7) (2004)
  689--703.
\newblock \href {http://dx.doi.org/10.1061/(ASCE)0733-9429(2004)130:7(689)}
  {\path{doi:10.1061/(ASCE)0733-9429(2004)130:7(689)}}.

\bibitem{zhao_etal_2016}
J.~Zhao, I.~Özgen Xian, R.~Hinkelmann, F.~Simons, D.~Liang, Comparison of
  capacity and non-capacity sediment transport models for dam break flow over
  movable bed, CRC Press, London, 2016, pp. 522--527.
\newblock \href {http://dx.doi.org/10.1201/9781315623207-96}
  {\path{doi:10.1201/9781315623207-96}}.

\bibitem{zhao_etal_2019}
J.~Zhao, I.~Özgen Xian, D.~Liang, T.~Wang, R.~Hinkelmann, A depth-averaged
  non-cohesive sediment transport model with improved discretization of flux
  and source terms, Journal of Hydrology 570 (2019) 647 -- 665.
\newblock \href {http://dx.doi.org/10.1016/j.jhydrol.2018.12.059}
  {\path{doi:10.1016/j.jhydrol.2018.12.059}}.

\bibitem{hu_etal_2019}
P.~Hu, Y.~Lei, J.~Han, Z.~Cao, H.~Liu, Z.~He, Computationally efficient
  modeling of hydro-sediment-morphodynamic processes using a hybrid local time
  step/global maximum time step, Advances in Water Resources 127 (2019) 26 --
  38.
\newblock \href {http://dx.doi.org/10.1016/j.advwatres.2019.03.006}
  {\path{doi:10.1016/j.advwatres.2019.03.006}}.

\bibitem{li_and_duffy_2011}
S.~Li, C.~J. Duffy, Fully coupled approach to modeling shallow water flow,
  sediment transport, and bed evolution in rivers, Water Resources Research
  47~(3).
\newblock \href {http://dx.doi.org/10.1029/2010WR009751}
  {\path{doi:10.1029/2010WR009751}}.

\bibitem{benkhaldoun_etal_2013}
F.~Benkhaldoun, I.~Elmahi, S.~Sari, M.~Seaid, An unstructured finite-volume
  method for coupled models of suspended sediment and bed load transport in
  shallow-water flows, International Journal for Numerical Methods in Fluids
  72~(9) (2013) 967--993.
\newblock \href {http://dx.doi.org/10.1002/fld.3771}
  {\path{doi:10.1002/fld.3771}}.

\bibitem{liu_etal_2015_0}
X.~Liu, J.~A.~I. Sedano, A.~Mohammadian, A robust coupled 2-{D} model for
  rapidly varying flows over erodible bed using central-upwind method with
  wetting and drying, Canadian Journal of Civil Engineering 42~(8) (2015)
  530--543.
\newblock \href {http://dx.doi.org/10.1139/cjce-2014-0524}
  {\path{doi:10.1139/cjce-2014-0524}}.

\bibitem{liu_etal_2015_1}
X.~Liu, A.~Mohammadian, A.~Kurganov, J.~A. {Infante Sedano}, Well-balanced
  central-upwind scheme for a fully coupled shallow water system modeling flows
  over erodible bed, Journal of Computational Physics 300 (2015) 202 -- 218.
\newblock \href {http://dx.doi.org/10.1016/j.jcp.2015.07.043}
  {\path{doi:10.1016/j.jcp.2015.07.043}}.

\bibitem{liu_and_beljadid_2017}
X.~Liu, A.~Beljadid, A coupled numerical model for water flow, sediment
  transport and bed erosion, Computers \& Fluids 154 (2017) 273 -- 284.
\newblock \href {http://dx.doi.org/10.1016/j.compfluid.2017.06.013}
  {\path{doi:10.1016/j.compfluid.2017.06.013}}.

\bibitem{xia_etal_2017}
C.~Xia, Z.~Cao, G.~Pender, A.~Borthwick, Numerical algorithms for solving
  shallow water hydro-sediment-morphodynamic equations, Engineering
  Computations 34 (2017) 00--00.
\newblock \href {http://dx.doi.org/10.1108/EC-01-2016-0026}
  {\path{doi:10.1108/EC-01-2016-0026}}.

\bibitem{kesserwani_etal_2014}
G.~Kesserwani, A.~Shamkhalchian, M.~J. Zadeh, {Fully Coupled Discontinuous
  Galerkin Modeling of Dam-Break Flows over Movable Bed with Sediment
  Transport}, Journal of Hydraulic Engineering 140~(4) (2014) 06014006.
\newblock \href {http://dx.doi.org/10.1061/(ASCE)HY.1943-7900.0000860}
  {\path{doi:10.1061/(ASCE)HY.1943-7900.0000860}}.

\bibitem{clare_etal_2020}
M.~Clare, J.~Percival, A.~Angeloudis, C.~Cotter, M.~Piggott,
  Hydro-morphodynamics {2D} modelling using a discontinuous {G}alerkin
  discretisation (Jan 2020).
\newblock \href {http://dx.doi.org/10.31223/osf.io/tpqvy}
  {\path{doi:10.31223/osf.io/tpqvy}}.

\bibitem{zhao_etal_1994}
D.~H. Zhao, H.~W. Shen, G.~Q. Tabios, J.~S. Lai, W.~Y. Tan, {Finite-Volume
  Two-Dimensional Unsteady-Flow Model for River Basins}, Journal of Hydraulic
  Engineering 120~(7) (1994) 863--883.
\newblock \href {http://dx.doi.org/10.1061/(ASCE)0733-9429(1994)120:7(863)}
  {\path{doi:10.1061/(ASCE)0733-9429(1994)120:7(863)}}.

\bibitem{anastasiou_and_chan_1997}
K.~Anastasiou, C.~T. Chan, {Solution of the 2D shallow water equations using
  the finite volume method on unstructured triangular meshes}, International
  Journal for Numerical Methods in Fluids 24~(11) (1997) 1225--1245.
\newblock \href
  {http://dx.doi.org/10.1002/(SICI)1097-0363(19970615)24:11<1225::AID-FLD540>3.0.CO;2-D}
  {\path{doi:10.1002/(SICI)1097-0363(19970615)24:11<1225::AID-FLD540>3.0.CO;2-D}}.

\bibitem{sleigh_etal_1998}
P.~Sleigh, P.~Gaskell, M.~Berzins, N.~Wright, An unstructured finite-volume
  algorithm for predicting flow in rivers and estuaries, Computers \& Fluids
  27~(4) (1998) 479 -- 508.
\newblock \href {http://dx.doi.org/10.1016/S0045-7930(97)00071-6}
  {\path{doi:10.1016/S0045-7930(97)00071-6}}.

\bibitem{aizinger_and_dawson_2002}
V.~Aizinger, C.~Dawson, A discontinuous {G}alerkin method for two-dimensional
  flow and transport in shallow water, Advances in Water Resources 25~(1)
  (2002) 67 -- 84.
\newblock \href {http://dx.doi.org/10.1016/S0309-1708(01)00019-7}
  {\path{doi:10.1016/S0309-1708(01)00019-7}}.

\bibitem{yoon_and_kang_2004}
T.~H. Yoon, S.-K. Kang, {Finite Volume Model for Two-Dimensional Shallow Water
  Flows on Unstructured Grids}, Journal of Hydraulic Engineering 130~(7) (2004)
  678--688.
\newblock \href {http://dx.doi.org/10.1061/(ASCE)0733-9429(2004)130:7(678)}
  {\path{doi:10.1061/(ASCE)0733-9429(2004)130:7(678)}}.

\bibitem{kubatko_etal_2006_nswe}
E.~J. Kubatko, J.~J. Westerink, C.~Dawson, \textit{hp} {D}iscontinuous
  {G}alerkin methods for advection dominated problems in shallow water flow,
  Computer Methods in Applied Mechanics and Engineering 196~(1) (2006) 437 --
  451.
\newblock \href {http://dx.doi.org/10.1016/j.cma.2006.05.002}
  {\path{doi:10.1016/j.cma.2006.05.002}}.

\bibitem{samii_etal_2019}
A.~Samii, K.~Kazhyken, C.~Michoski, C.~Dawson, A {C}omparison of the {E}xplicit
  and {I}mplicit {H}ybridizable {D}iscontinuous {G}alerkin {M}ethods for
  {N}onlinear {S}hallow {W}ater {E}quations, Journal of Scientific Computing
  80~(3) (2019) 1936--1956.
\newblock \href {http://dx.doi.org/10.1007/s10915-019-01007-z}
  {\path{doi:10.1007/s10915-019-01007-z}}.

\bibitem{bremer_etal_2019}
M.~Bremer, K.~Kazhyken, H.~Kaiser, C.~Michoski, C.~Dawson, Performance
  {C}omparison of {HPX} {V}ersus {T}raditional {P}arallelization {S}trategies
  for the {D}iscontinuous {G}alerkin {M}ethod, Journal of Scientific Computing
  80~(2) (2019) 878--902.
\newblock \href {http://dx.doi.org/10.1007/s10915-019-00960-z}
  {\path{doi:10.1007/s10915-019-00960-z}}.

\bibitem{green_naghdi_1976}
A.~E. Green, P.~M. Naghdi, A derivation of equations for wave propagation in
  water of variable depth, Journal of Fluid Mechanics 78~(2) (1976) 237–246.
\newblock \href {http://dx.doi.org/10.1017/S0022112076002425}
  {\path{doi:10.1017/S0022112076002425}}.

\bibitem{chazel_etal_2011}
F.~Chazel, D.~Lannes, F.~Marche, {Numerical Simulation of Strongly Nonlinear
  and Dispersive Waves Using a Green-Naghdi Model}, Journal of Scientific
  Computing 48~(1) (2011) 105--116.
\newblock \href {http://dx.doi.org/10.1007/s10915-010-9395-9}
  {\path{doi:10.1007/s10915-010-9395-9}}.

\bibitem{bonneton_etal_2011}
P.~Bonneton, F.~Chazel, D.~Lannes, F.~Marche, M.~Tissier, A splitting approach
  for the fully nonlinear and weakly dispersive {G}reen–{N}aghdi model,
  Journal of Computational Physics 230~(4) (2011) 1479 -- 1498.
\newblock \href {http://dx.doi.org/10.1016/j.jcp.2010.11.015}
  {\path{doi:10.1016/j.jcp.2010.11.015}}.

\bibitem{panda_etal_2014}
N.~Panda, C.~Dawson, Y.~Zhang, A.~B. Kennedy, J.~J. Westerink, A.~S. Donahue,
  {Discontinuous Galerkin methods for solving Boussinesq–Green–Naghdi
  equations in resolving non-linear and dispersive surface water waves},
  Journal of Computational Physics 273 (2014) 572 -- 588.
\newblock \href {http://dx.doi.org/10.1016/j.jcp.2014.05.035}
  {\path{doi:10.1016/j.jcp.2014.05.035}}.

\bibitem{lannes_and_marche_2015}
D.~Lannes, F.~Marche, A new class of fully nonlinear and weakly dispersive
  {G}reen–{N}aghdi models for efficient 2{D} simulations, Journal of
  Computational Physics 282 (2015) 238 -- 268.
\newblock \href {http://dx.doi.org/10.1016/j.jcp.2014.11.016}
  {\path{doi:10.1016/j.jcp.2014.11.016}}.

\bibitem{duran_marche_2015}
A.~Duran, F.~Marche, {Discontinuous-Galerkin Discretization of a New Class of
  Green-Naghdi Equations}, Communications in Computational Physics 17~(3)
  (2015) 721–760.
\newblock \href {http://dx.doi.org/10.4208/cicp.150414.101014a}
  {\path{doi:10.4208/cicp.150414.101014a}}.

\bibitem{duran_and_marche_2017}
A.~Duran, F.~Marche, A discontinuous {G}alerkin method for a new class of
  {G}reen–{N}aghdi equations on simplicial unstructured meshes, Applied
  Mathematical Modelling 45 (2017) 840 -- 864.
\newblock \href {http://dx.doi.org/10.1016/j.apm.2017.01.030}
  {\path{doi:10.1016/j.apm.2017.01.030}}.

\bibitem{samii_and_dawson_2018}
A.~Samii, C.~Dawson, An explicit hybridized discontinuous {G}alerkin method for
  {S}erre–{G}reen–{N}aghdi wave model, Computer Methods in Applied
  Mechanics and Engineering 330 (2018) 447 -- 470.
\newblock \href {http://dx.doi.org/10.1016/j.cma.2017.11.001}
  {\path{doi:10.1016/j.cma.2017.11.001}}.

\bibitem{marche_2020}
F.~Marche, {Combined Hybridizable Discontinuous Galerkin (HDG) and Runge-Kutta
  Discontinuous Galerkin (RK-DG) formulations for Green-Naghdi equations on
  unstructured meshes}, Journal of Computational Physics 418 (2020) 109637.
\newblock \href {http://dx.doi.org/10.1016/j.jcp.2020.109637}
  {\path{doi:10.1016/j.jcp.2020.109637}}.

\bibitem{li_duffy_2011}
S.~Li, C.~J. Duffy, Fully coupled approach to modeling shallow water flow,
  sediment transport, and bed evolution in rivers, Water Resources Research
  47~(3).
\newblock \href {http://dx.doi.org/10.1029/2010WR009751}
  {\path{doi:10.1029/2010WR009751}}.

\bibitem{diaz_etal_2008}
{M. J. Castro Díaz}, {E. D. Fernández-Nieto}, {A. M. Ferreiro}, Sediment
  transport models in {S}hallow {W}ater equations and numerical approach by
  high order finite volume methods, Computers {\&} Fluids 37~(3) (2008) 299 --
  316.
\newblock \href {http://dx.doi.org/10.1016/j.compfluid.2007.07.017}
  {\path{doi:10.1016/j.compfluid.2007.07.017}}.

\bibitem{cordier_etal_2011}
S.~Cordier, M.~Le, {T. Morales de Luna}, Bedload transport in shallow water
  models: {W}hy splitting (may) fail, how hyperbolicity (can) help, Advances in
  Water Resources 34~(8) (2011) 980 -- 989.
\newblock \href {http://dx.doi.org/10.1016/j.advwatres.2011.05.002}
  {\path{doi:10.1016/j.advwatres.2011.05.002}}.

\bibitem{grass_1981}
{A. J. Grass}, Sediment {T}ransport by {W}aves and {C}urrents, SERC London
  Centre for Marine Technology, Report No. FL29.

\bibitem{strang_1968}
G.~Strang, On the {C}onstruction and {C}omparison of {D}ifference {S}chemes,
  SIAM Journal on Numerical Analysis 5~(3) (1968) 506--517.
\newblock \href {http://dx.doi.org/10.1137/0705041}
  {\path{doi:10.1137/0705041}}.

\bibitem{mirabito_etal_2011}
C.~Mirabito, C.~Dawson, {E. J. Kubatko}, {J. J. Westerink}, S.~Bunya,
  Implementation of a discontinuous {G}alerkin morphological model on
  two-dimensional unstructured meshes, Computer Methods in Applied Mechanics
  and Engineering 200~(1) (2011) 189 -- 207.
\newblock \href {http://dx.doi.org/10.1016/j.cma.2010.08.004}
  {\path{doi:10.1016/j.cma.2010.08.004}}.

\bibitem{harten_etal_1983}
A.~Harten, P.~D. Lax, B.~v. Leer, On {U}pstream {D}ifferencing and
  {G}odunov-{T}ype {S}chemes for {H}yperbolic {C}onservation {L}aws, SIAM
  Review 25~(1) (1983) 35--61.
\newblock \href {http://dx.doi.org/10.1137/1025002}
  {\path{doi:10.1137/1025002}}.

\bibitem{nguyen_peraire_2012}
N.~Nguyen, J.~Peraire, Hybridizable discontinuous {G}alerkin methods for
  partial differential equations in continuum mechanics, Journal of
  Computational Physics 231~(18) (2012) 5955 -- 5988.
\newblock \href {http://dx.doi.org/10.1016/j.jcp.2012.02.033}
  {\path{doi:10.1016/j.jcp.2012.02.033}}.

\bibitem{bunya_etal_2009}
S.~Bunya, E.~J. Kubatko, J.~J. Westerink, C.~Dawson, A wetting and drying
  treatment for the {R}unge–{K}utta discontinuous {G}alerkin solution to the
  shallow water equations, Computer Methods in Applied Mechanics and
  Engineering 198~(17) (2009) 1548 -- 1562.
\newblock \href {http://dx.doi.org/10.1016/j.cma.2009.01.008}
  {\path{doi:10.1016/j.cma.2009.01.008}}.

\bibitem{krivodonova_etal_2004}
L.~Krivodonova, J.~Xin, J.-F. Remacle, N.~Chevaugeon, J.~Flaherty, Shock
  detection and limiting with discontinuous {G}alerkin methods for hyperbolic
  conservation laws, Applied Numerical Mathematics 48~(3) (2004) 323 -- 338.
\newblock \href {http://dx.doi.org/10.1016/j.apnum.2003.11.002}
  {\path{doi:10.1016/j.apnum.2003.11.002}}.

\bibitem{cockburn_shu_2001}
B.~Cockburn, C.~Shu, Runge-{K}utta discontinuous {G}alerkin methods for
  convection-dominated problems, Journal of Scientific Computing 16~(3) (2001)
  173--261.
\newblock \href {http://dx.doi.org/10.1023/A:1012873910884}
  {\path{doi:10.1023/A:1012873910884}}.

\bibitem{eigen}
G.~Guennebaud, B.~Jacob, et~al., Eigen v3, \url{www.eigen.tuxfamily.org}
  (2010).

\bibitem{blaze}
K.~{Iglberger}, Blaze {C}++ {L}inear {A}lgebra {L}ibrary,
  \url{www.bitbucket.org/blaze-lib} (2012).

\bibitem{petsc}
S.~Balay, S.~Abhyankar, M.~F. Adams, J.~Brown, P.~Brune, K.~Buschelman,
  L.~Dalcin, A.~Dener, V.~Eijkhout, W.~D. Gropp, D.~Karpeyev, D.~Kaushik, M.~G.
  Knepley, D.~A. May, L.~C. McInnes, R.~T. Mills, T.~Munson, K.~Rupp, P.~Sanan,
  B.~F. Smith, S.~Zampini, H.~Zhang, H.~Zhang, Portable, {E}xtensible {T}oolkit
  for {S}cientific {C}omputation, \url{www.mcs.anl.gov/petsc} (2019).

\bibitem{hpx}
H.~Kaiser, B.~A. Lelbach, T.~Heller, M.~Simberg, A.~Bergé, J.~Biddiscombe,
  A.~Bikineev, G.~Mercer, A.~Schäfer, K.~Huck, A.~S. Lemoine, T.~Kwon,
  J.~Habraken, M.~Anderson, M.~Copik, S.~R. Brandt, M.~Stumpf, D.~Bourgeois,
  D.~Blank, S.~Jakobovits, V.~Amatya, L.~Viklund, Z.~Khatami, P.~Diehl,
  T.~Pathak, D.~Bacharwar, S.~Yang, E.~Schnetter, {STEllAR-GROUP/hpx: HPX
  V1.4.1: The C++ Standards Library for Parallelism and Concurrency} (Feb.
  2020).
\newblock \href {http://dx.doi.org/10.5281/zenodo.3675272}
  {\path{doi:10.5281/zenodo.3675272}}.

\bibitem{dubiner_1991}
M.~Dubiner, Spectral methods on triangles and other domains, Journal of
  Scientific Computing 6~(4) (1991) 345–390.
\newblock \href {http://dx.doi.org/10.1007/BF01060030}
  {\path{doi:10.1007/BF01060030}}.

\bibitem{fraccarollo_capart_2002}
L.~Fraccarollo, H.~Capart, Riemann wave description of erosional dam-break
  flows, Journal of Fluid Mechanics 461 (2002) 183–228.
\newblock \href {http://dx.doi.org/10.1017/S0022112002008455}
  {\path{doi:10.1017/S0022112002008455}}.

\bibitem{capart_young_1998}
H.~Capart, D.~L. Young, Formation of a jump by the dam-break wave over a
  granular bed, Journal of Fluid Mechanics 372 (1998) 165–187.
\newblock \href {http://dx.doi.org/10.1017/S0022112098002250}
  {\path{doi:10.1017/S0022112098002250}}.

\bibitem{goutiere_etal_2011}
L.~Goutiere, S.~{Soares-Fraz\~ao}, Y.~Zech, Dam-break flow on mobile bed in
  abruptly widening channel: experimental data, Journal of Hydraulic Research
  49~(3) (2011) 367--371.
\newblock \href {http://dx.doi.org/10.1080/00221686.2010.548969}
  {\path{doi:10.1080/00221686.2010.548969}}.

\bibitem{soares_etal_2012}
S.~{Soares-Fraz\~ao}, R.~Canelas, Z.~Cao, L.~Cea, H.~M. Chaudhry, A.~{Die
  Moran}, K.~{El Kadi}, R.~Ferreira, I.~{Fraga Cadórniga},
  N.~Gonzalez-Ramirez, M.~Greco, W.~Huang, J.~Imran, J.~{Le Coz}, R.~Marsooli,
  A.~Paquier, G.~Pender, M.~Pontillo, J.~Puertas, B.~Spinewine,
  C.~Swartenbroekx, R.~Tsubaki, C.~Villaret, W.~Wu, Z.~Yue, Y.~Zech, Dam-break
  flows over mobile beds: experiments and benchmark tests for numerical models,
  Journal of Hydraulic Research 50~(4) (2012) 364--375.
\newblock \href {http://dx.doi.org/10.1080/00221686.2012.689682}
  {\path{doi:10.1080/00221686.2012.689682}}.

\bibitem{sumer_etal_2011}
B.~M. Sumer, M.~B. Sen, I.~Karagali, B.~Ceren, J.~Fredsøe, M.~Sottile,
  L.~Zilioli, D.~R. Fuhrman, Flow and sediment transport induced by a plunging
  solitary wave, Journal of Geophysical Research: Oceans 116~(C1).
\newblock \href {http://dx.doi.org/10.1029/2010JC006435}
  {\path{doi:10.1029/2010JC006435}}.

\bibitem{young_etal_2010}
Y.~L. Young, H.~Xiao, T.~Maddux, Hydro- and morpho-dynamic modeling of breaking
  solitary waves over a fine sand beach. {P}art {I}: {E}xperimental study,
  Marine Geology 269~(3) (2010) 107 -- 118.
\newblock \href {http://dx.doi.org/10.1016/j.margeo.2009.12.009}
  {\path{doi:10.1016/j.margeo.2009.12.009}}.

\end{thebibliography}

\end{document}